\documentclass
[
a4paper,
pdftex,
onecolumn,
DIV=11,
abstracton
]
{scrartcl}

\usepackage
{
	graphicx,
	amssymb,
	amsmath,
	amsthm,
	xcolor,
	dsfont,
	algpseudocode,
	authblk,
	enumerate,
	bm,
	mathtools,
	float
}
\usepackage[pdffitwindow=false,
            plainpages=false,
            pdfpagelabels=true,
            pdfpagemode=UseOutlines,
            pdfpagelayout=SinglePage,
            bookmarks=false,
            colorlinks=true,
            hyperfootnotes=false,
            linkcolor=blue,
            urlcolor=blue!30!black,
            citecolor=green!50!black]{hyperref}

\usepackage[left=3.5cm,right=3.5cm,top=2cm]{geometry}

\definecolor{upbgray}{rgb}{0.780392156863, 0.788235294118, 0.780392156863}
\definecolor{upbblue}{rgb}{0, 0.125490196078, 0.3568627451}
\definecolor{upbblue2}{rgb}{0, 0.623529412, 0.874509804}

\usepackage{algorithm}

\usepackage[bf]{caption}

\usepackage{subfigure}

\DeclareMathAlphabet{\mathpzc}{OT1}{pzc}{m}{it}

\newcommand{\R}{\mathbb{R}}

\newcommand{\N}{\mathbb{N}}
\newcommand{\cB}{{\mathcal{B}}}
\newcommand{\cC}{{\mathcal{C}}}
\newcommand{\cM}{{\mathcal{M}}}

\newcommand{\cP}{{\mathcal{P}}}
\newcommand{\cA}{{\mathcal{A}}}
\newcommand{\cL}{{\mathcal{L}}}
\newcommand{\cW}{{\mathcal{W}}}

\newcommand{\setdist}[1]{\textrm{dist}\left(#1\right)}

\providecommand{\norm}[1]{\left\lVert #1 \right\rVert}

\DeclareMathOperator{\diam}{diam}

\DeclareMathOperator{\vol}{vol}

\newtheorem{theorem}{Theorem}[section]

\newtheorem{proposition}[theorem]{Proposition}

\newtheorem{remark}[theorem]{Remark}

\makeatletter
\renewcommand*\env@matrix[1][*\c@MaxMatrixCols c]{%
	\hskip -\arraycolsep
	\let\@ifnextchar\new@ifnextchar
	\array{#1}}
\makeatother

\newcommand{\old}[1]{\textcolor{green}{{ }}}

\begin{document}
\title{Revealing the intrinsic geometry of finite dimensional invariant sets of infinite dimensional dynamical systems}

\author[1]{Raphael Gerlach}
\author[2]{P\'eter Koltai}
\author[1]{Michael Dellnitz}
\affil[1]{\normalsize Department of Mathematics, Paderborn University, 33098 Paderborn, Germany.}
\affil[2]{\normalsize Institute of Mathematics, Freie Universit\"at Berlin, 14195 Berlin, Germany.}
	
	\maketitle
%
\begin{abstract}\label{sec:abstract}
	Embedding techniques allow the approximations of finite dimensional attractors and manifolds of infinite dimensional dynamical systems via subdivision 
	and continuation
	methods. These approximations give a topological one-to-one image of the original set. In order to additionally reveal their geometry we use diffusion maps
	to find intrinsic coordinates. We illustrate our results on the unstable manifold of the one-dimensional Kuramoto--Sivashinsky equation, as well as for the attractor of the Mackey--Glass delay differential equation.
\end{abstract}
\section{Introduction}\label{sec:introduction}
For the understanding of the long term behavior of non-linear complicated dynamical systems the objects of interest are invariant sets such as the global attractor and invariant manifolds. To numerically approximate these sets set-oriented methods have been developed~\cite{DH96, DH97, DJ99, FD03}. The underlying idea is to cover the set by outer approximations that are generated by multilevel subdivision or continuation methods. They have been used successfully in various application areas such as molecular dynamics~\cite{SHD01}, astrodynamics~\cite{DJLMPPRT05} or ocean dynamics~\cite{FHRSSG12}.

Recently these methods have been extended from the finite dimensional setting to the treatment of infinite dimensional dynamical systems such as partial differential equations. In particular, in~\cite{DHZ16} the subdivision algorithm developed in~\cite{DH97} was adapted to allow the approximation of finite dimensional attractors of semi-flows in (possible infinite dimensional) Banach spaces. Additionally, finite dimensional manifolds of steady states can be computed by the extension of the continuation algorithm introduced in~\cite{DH96} to infinite dimensional systems~\cite{ZDG18}. Both of the methods rely on embbedding techniques~\cite{whitney1936differentiable,takens1981detecting,HuntKaloshin99,R05} that allow the construction of the so--called \emph{core dynamical system}, that is a finite dimensional system which is topologically conjugate to the original dynamics on the attractor. Thus, the traditional set-oriented numerical methods can be used to approximate \emph{embedded} attractors or \emph{embedded} manifolds which are one-to-one images of the corresponding set in the Banach space.

Beyond these topological objects, the finite-dimensional core dynamical system is also suitable for calculating \emph{measure-theoretic dynamical} quantities associated with the dynamics, like invariant measures~\cite{DHZ16}. Other dynamical characteristics, such as (quasi-)periodic and (almost-)invariant motion can also be analyzed~\cite{DJ99}. This is also the goal of data-driven Koopman-operator approaches~\cite{giannakis2018koopman,arbabi2017study} and Dynamic Mode Decomposition~\cite{SchSe08,WKR15}---without delivering any claims about the topology of the attractor.

Driven by the desire to obtain further intuitive understanding of the geometric structure of the attractors (and thus hopefully also of the dynamics on them) of infinite dimensional systems, here we will identify nonlinear coordinates revealing their intrinsic geometry in the embedding space.
To this end, we are using diffusion maps to obtain the geometric and dynamical structure of the covering of an embedded attractor or manifold. Diffusion maps are one among many data--driven manifold learning techniques that find intrinsic coordinates of a data set~\cite{TdSL00,RoSa00,DoGr03,BeNi03,ZhZh04}. First introduced by Coifman and Lafon \cite{coifman2006diffusion, coifman2005MM, coifman2005DMAPS}, diffusion maps is a nonlinear feature extraction algorithm that computes a family of embeddings of a (possibly) high dimensional data set into a low dimensional space, whose coordinates are given by the eigenvectors and eigenvalues of a diffusion operator on the data. Different from linear dimensionality reductions methods, such as principal component analysis (POD), diffusion maps focus on discovering the underlying manifold from which the data set is sampled. In addition to that the algorithm is robust to noise perturbation such that it can deal with the outer approximations that cover the set of interest generated by the set-oriented numerical methods.

A detailed outline of this paper is as follows. In Section \ref{sec:SubdivisonContinuation} we briefly summarize the results of \cite{DHZ16} and \cite{ZDG18}. We state the main embedding results of \cite{HuntKaloshin99} and \cite{R05} and we describe the construction of the \emph{core dynamical system} on the observation space. Afterwards, we explain how the classical subdivision and continuation algorithm is extended to the infinite dimensional setting. In Section \ref{sec:Diffusion Maps} we review the concept of diffusion maps and give a suitable numerical implementation for our purposes. Finally, in Section~\ref{sec:applications} we apply this method to the embedded unstable manifold for a one dimensional Kuramoto--Sivashinksy equation for different parameter values and the embedded attractor of the Mackey--Glass delay differential equation.

\section{Review of Subdivision and Continuation methods for infinite dimensional dynamical systems}\label{sec:SubdivisonContinuation}

Since we want to analyze the geometry of invariant sets of infinite dimensional dynamical systems we start with a short review of the novel set-oriented methods to approximate those sets.

\subsection{Embedding Techniques}\label{ssec:embeddings}
We consider dynamical systems of the form
\begin{equation}\label{eq:DS}
u_{j+1} = \Phi (u_j),\quad j=0,1,\ldots,
\end{equation}
where $\Phi:Y\rightarrow Y$ is Lipschitz continuous on a Banach space $Y$.
Moreover, we assume that $\Phi$ has an invariant compact set $\cA$, that is
\[
\Phi(\cA) = \cA.
\]

In order to approximate invariant subsets of $\cA$ or $\cA$ itself  we combine classical subdivision and continuation techniques for the computation of such objects in a finite dimensional space with infinite dimensional embedding results (cf. \cite{HuntKaloshin99,R05}). The following theorems allow us to map $\cA$ into a finite dimensional space $\R^k$ such that this map is generically---in the sense of \emph{prevalence}\footnote{A Borel subset $S$ of a normed linear space $V$ is	\emph{prevalent} if there is a finite dimensional subspace $E$ of $V$ (the `probe space') such that for each $v \in V,\ v+e$ belongs to $S$ for (Lebesgue) almost every $e\in E$.}
\cite{SYC91}---one-to-one on $\cA$. To do so the embedding dimension $k\in\N$ has to be chosen large enough depending on the \emph{upper box counting dimension} $d_{\mathrm{box}}$ and \emph{thickness exponent}\footnote{The \emph{thickness exponent} measures roughly speaking, how well $\cA$ can be approximated be finite dimensional linear subspaces}~$\sigma$~\cite{HuntKaloshin99}.

\begin{theorem}[\cite{HuntKaloshin99}]
	Let $Y$ be a Banach space and $\cA\subset Y$ compact, with upper box counting dimension $d = d_{\mathrm{box}}(\cA; Y)$ and thickness exponent $\sigma = \sigma(\cA; Y)$. Let $N>2d$ be an integer, and let $\alpha \in\R$ with
	\[ 0 < \alpha < \frac{ N-2d }{ N\cdot(1+\sigma) }. \]
	Then, for a prevalent set of bounded linear maps $\cL:Y\to\R^N$ there is $C>0$ such that
	\[ C\cdot \| \cL(x-y)\|^\alpha \geq \|x-y\| \quad\text{for all }x,
	y\in \cA. \]
	\label{thm:HK99}
\end{theorem}

This theorem lays the foundation for Robinson's main result concerning delay embedding techniques.

\begin{theorem}[\cite{R05}]\label{thm:R05}
	Let $Y$ be a Banach space and $\cA \subset Y$ a compact, invariant set,
	with upper box counting dimension $d$, and thickness exponent $\sigma$. Choose an integer $k> 2(1+\sigma)d$ and suppose further that the set $A_p$
	of $p$-periodic points of $\Phi$ satisfies $d_{\mathrm{box}}(A_p;Y) < p/(2+2\sigma)$ for $p=1,\ldots,k$.
	Then, for a prevalent set of Lipschitz maps $f:Y\rightarrow \R$ the observation map
	$D_k[f,\Phi] : Y \rightarrow \R^k$ defined by
	\begin{equation*}\label{def:D_k}
	D_k[f,\Phi] (u) = (f(u),f(\Phi(u)),\ldots,f(\Phi^{k-1}(u)))^T
	\end{equation*}
	is one-to-one on $\cA$.
\end{theorem}

This result can be generalized to the case where several different observables are evaluated. In fact, we can use $k\in \mathbb{N}$ observables $f_i:Y\to \R$, $i=1,\ldots,k$ such that the \emph{observation map} $R:Y\to \R^k$ given by
\begin{align}\label{eq:R}
R(u) = R[f](u) = (f_1(u),\ldots,f_k(u))^T.
\end{align}
is one-to-one on~$\mathcal{A}$. Moreover, it is reasonable to assume that the thickness exponent is zero \cite{friz1999smooth}. Thus, provided $k>2d$, the observation map $R$ is generically---in the sense of prevalence---one-to-one on~$\mathcal{A}$. With that in mind, we will compute the image $R(\mathcal{A})$ of $\cA$ under the observation map instead of $\mathcal{A}$ itself. Likewise we aim to approximate embedded manifolds $R(\mathcal{W}_u(u^*))\subset R(\mathcal{A})$ for some unstable steady state~$u^*\in \cA$.

\subsection{The core dynamical system (CDS)}
\label{sec:comp_emb_att}

Using the results from Section~\ref{ssec:embeddings} a finite dimensional dynamical system $\varphi$, the so--called \emph{core dynamical system} (CDS), can be created that essentially has the same dynamics as the infinite dimensional system~\eqref{eq:DS}. In this section we will briefly review the construction of the CDS.

\begin{figure}[H]
	\centering
	\includegraphics[width=0.56\textwidth]{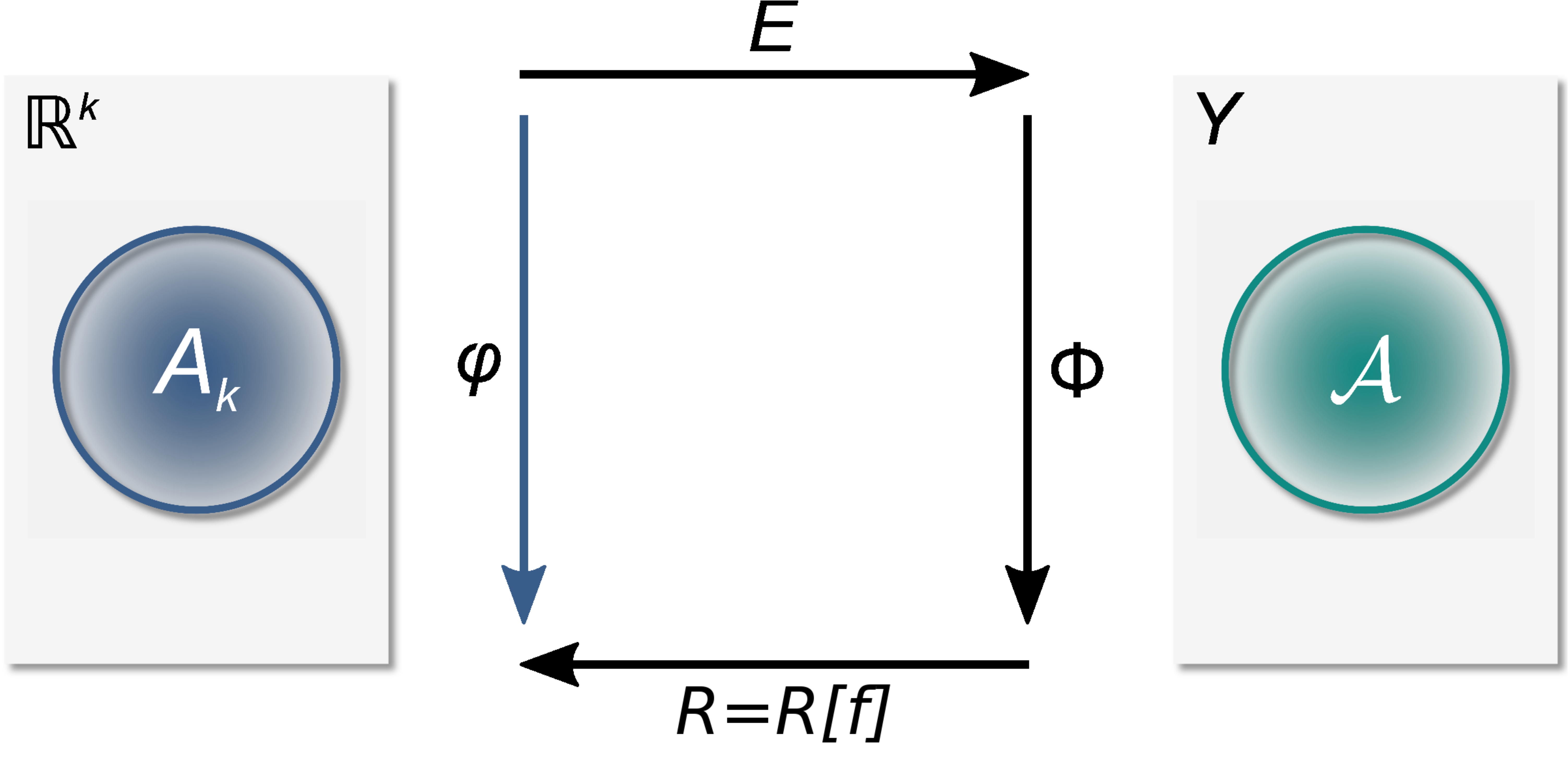} 
	\caption{Definition of the core dynamical system~$\varphi$.}
	\label{fig:kommdia}
\end{figure}

Denote by $A_k$ the image of $\cA\subset Y$ under the observation map $R:Y \rightarrow \R^k$, that is
\begin{align*}
A_k = R(\cA),
\end{align*}
where $R$ is defined in~\eqref{eq:R}. The CDS is then constructed as follows:
First we define $\varphi$ on the set $A_k$ by
\begin{equation*}
\label{def:phi}
\varphi = R \circ \Phi \circ \widetilde E,
\end{equation*}
where $\widetilde E:A_k \rightarrow Y$ is the continuous map satisfying
\begin{equation*}
\label{eq:condRE}
(\widetilde E\circ R)(u) = u\quad \forall u\in \cA
\quad\text{ and }\quad
(R\circ \widetilde E)(x) = x\quad \forall x\in A_k.
\end{equation*}
This is possible due to the fact that $R$ is invertible as a mapping from $\cA$ to $A_k$. Using a generalization of Tietze's extension theorem \cite{Dugundji51} we extend $\widetilde{E}$ to a continuous map $E:\R^k\to Y$ with $E\vert_{A_k} = \widetilde E$ (see Figure~\ref{fig:kommdia}) bringing us in the position to define the CDS $\varphi$ on $\R^k$, i.e., 
\begin{equation*}\label{eq:CDS}
x_{j+1} = \varphi(x_j), \quad j=0,1,2,\ldots,
\end{equation*}
where $\varphi = R \circ \Phi \circ E: \R^k \rightarrow \R^k$. Note that by construction the dynamics of the CDS $\varphi$ on $A_k$ is topologically conjugate to that of $\Phi$ on~$\cA$.
\begin{proposition}[{\cite[Proposition 1]{DHZ16}}]
	\label{prop:phi_is_cont}
	There is a continuous map ${\varphi:\R^k \to \R^k}$ satisfying
	\begin{equation*}
	\varphi(R(u)) = R(\Phi(u)) \text{ for all } u\in \cA.
	\label{eq:conjugonA}
	\end{equation*}
\end{proposition}

\subsection{Computation of embedded attractors via subdivision}
\label{ssec:subdiv}
Now we shall give a brief review of the adapted subdivision scheme developed in~\cite{DHZ16} that allows us to approximate the set~$A_k$.

Let $Q \subset \R^k$ be a compact set and suppose $A_k\subset Q$ for simplicity. The \emph{global attractor relative to} $Q$ is defined by
\begin{equation*}
\label{eq:relativeAttractor}
A_Q = \bigcap_{j\ge 0} \varphi^j(Q).
\end{equation*}
The aim is to approximate this set with a subdivision procedure. Given an initial finite collection $\cB_0$ of compact subsets of $\R^k$ such that
\[
Q=\bigcup_{B\in\cB_0} B,
\]
we recursively obtain $\cB_\ell$ from $\cB_{\ell-1}$
for $\ell=1,2,\ldots$ in two steps such that the diameter
\[
\diam(\cB_\ell) = \max_{B\in\cB_\ell}\diam(B)
\]
converges to zero for~$\ell\rightarrow\infty$.
\begin{algorithm}[H]
	\caption{The subdivision method for embedded global attractors}\label{alg:subdivision}
	\vspace{1em}
	\flushleft {\em Initialization:} Given $k>2(1+\sigma)d$ we choose a compact set $Q \subset \R^k$, such that $A_k \subset Q$. Fix $0<\theta_{\min} \le \theta_{\max} < 1$.
	\begin{enumerate}[1)]
		\item {\em Subdivision:} Construct a new collection
		$\hat\cB_\ell$ such that
		\begin{equation*}
		\label{eq:sd1}
		\bigcup_{B\in\hat\cB_\ell}B = \bigcup_{B\in\cB_{\ell-1}}B
		\end{equation*}
		and
		\begin{equation*}
		\label{eq:sd2}
		\diam(\hat\cB_\ell) = \theta_\ell\diam(\cB_{\ell-1}),
		\end{equation*}
		where $0<\theta_{\min} \le \theta_\ell\le \theta_{\max} < 1$.
		\item  {\em Selection:} Define the new collection $\cB_\ell$ by
		\begin{equation*}
		\label{eq:select}
		\cB_\ell=\left\{B\in\hat\cB_\ell : \exists \hat B\in\hat\cB_\ell
		~\mbox{such that}~\varphi^{-1}(B)\cap\hat B\ne\emptyset\right\}.
		\end{equation*}
	\end{enumerate}
	
\end{algorithm}

The first step is responsible for decreasing the size of the sets of increasing $\ell$. In fact, by construction
\begin{equation*}\label{eq:diamB}
\diam(\cB_\ell)\leq\theta_{\max}^\ell\diam(\cB_0)\rightarrow 0\quad
\mbox{for $\ell\rightarrow\infty$.}
\end{equation*}

In the second step each subset whose preimage does neither
intersect itself nor any other subset in $\hat\cB_\ell$ is removed. 
Denote by $Q_\ell$ the collection of compact subsets obtained after $\ell$
subdivision steps, that is 
\[
Q_\ell=\bigcup_{B\in\cB_\ell}B.
\]
Since the $Q_\ell$'s define a nested sequence of compact sets, that is, $Q_{\ell+1}\subset Q_\ell$ we conclude for each $m$
\begin{equation*}
Q_m = \bigcap\limits_{\ell =1}^{m} Q_\ell.
\end{equation*}
Then by considering 
\begin{equation*}
Q_{\infty} = \bigcap\limits_{\ell =1}^{\infty} Q_\ell
\end{equation*}
as the limit of the $Q_\ell$'s the selection step accounts for the fact that $Q_m$ approaches the relative global attractor.

\begin{proposition}[{\cite[Proposition 2]{DHZ16}}]\label{prop:AQ_Qinfty}
	Suppose that $A_Q$ satisfies ${\varphi^{-1}(A_Q) \subset A_Q}$. Then
	\[
	A_Q = Q_\infty.
	\]
\end{proposition}
We note that we can, in general, not expect that $A_k = A_Q$. In fact, by construction $A_Q$ may contain several invariant sets and related heteroclinic connections. However, if $\cA$ is an attracting set equality can be proven~(see \cite{DHZ16}).

\subsection{A subdivision and continuation technique for embedded unstable manifolds}\label{ssec:continuation}
In \cite{ZDG18} the classical continuation method of \cite{DH96} has been extended to the approximation of embedded unstable manifolds. In the following  we state the main result of this scheme. Let us denote by
\begin{equation*}\label{eq:unstable_manifold}
\cW_{\Phi}^u(u^*) \subset \cA
\end{equation*}
the unstable manifold of $u^* \in \cA$, where $u^*$ is a steady state solution of the infinite dimensional dynamical system $\Phi$ (cf. \eqref{eq:DS}). Furthermore, let us define the \emph{embedded unstable manifold} $W^u(p)$ by
\begin{equation*}\label{eq:embedded_unstable_manifold}
W^u(p) = R(\cW_{\Phi}^u(u^*)) \subset A_k,
\end{equation*}
where $p = R(u^*)$ and $R$ is the observation map introduced in Section~\ref{ssec:embeddings}. Choose a compact set $Q \subset \R^k$ containing $p$ and we assume for simplicity that $Q$ is large enough so that it contains the entire closure of the embedded unstable manifold, i.e.,
\begin{equation*}
\label{eq:closure}
\overline{W^u(p)} \subset Q.
\end{equation*}
For the purpose of initializing the developed algorithm we define a partition $\cP$ of $Q$ to be a finite family of compact subsets of $Q$ such that
\[
\bigcup_{B \in \cP}B = Q\quad \mbox{and}\quad \mbox{int} B\cap \mbox{int} B' = \emptyset,\ \mbox{for all } B, B'\in \cP,~ B\neq B'.
\]
Moreover, we denote by $\cP(x) \in \cP$ the element of $\cP$ containing $x\in Q$. We consider a nested sequence $\cP_s,\ s \in \N$, of successively finer partitions of $Q$, requiring that for all $B \in \cP_s$ there exist $B_1,\ldots,B_m \in \cP_{s+1}$ such that $B = \cup_i B_i$ and ${\diam(B_i)\leq\theta\diam(B)}$ for some $0 < \theta < 1$. A set $B \in \cP_s$ is said to be of \emph{level} $s$. 

The aim of the continuation method is to approximate subsets $W_j \subset W^u(p)$ where $W_0 = W_{loc}^u(p) =R(\overline{\cW_{\Phi,loc}^u(u^*)})$ is the local embedded unstable manifold and
\[
W_{j+1} = \varphi(W_j) \quad \mbox{for $j = 0,1,2,\ldots$}
\]
in two steps:

At first we use the subdivision Algorithm \ref{alg:subdivision} to approximate the local embedded unstable manifold $W_{loc}^u(p) = R(\overline{\cW_{\Phi,loc}^u(u^*)})$. To this end, we compute the relative global attractor $A_C$ of a compact neighborhood $C\subset A_k$. The idea of the continuation algorithm is then to globalize this local covering of $W_{loc}^u(p)$ to obtain an approximation of the compact subsets $W_j\subset W^u(p)$ or even the entire closure $\overline{W^u(p)}$.

\begin{algorithm}[H]
	\caption{The continuation method for embedded unstable manifolds}\label{alg:continuation}
	\vspace{1em}
	{\em Initialization:}  Given $k>2(1+\sigma)d$ we choose an initial box $Q \subset \R^k$, such that $A_k \subset Q$. Choose a partition $\cP_s$ of $Q$ and a set $C \in \cP_s$ such that $p = R(u^*) \in C$.
	\begin{enumerate}[1)]
		\item Apply the subdivision algorithm with $\ell$ subdivision steps to $\cB_0 = \{ C \}$ to obtain a covering $\cB_\ell \subset \cP_{s+\ell}$ of the local embedded unstable manifold $A_{C}$. 
		\item Set \
		\[
		\cC_0^{(\ell)} = \cB_\ell.
		\]
		\item For $j=0,1,2,\ldots$ define 
		\begin{align*}\label{eq:continuation}
		\cC_{j+1}^{(\ell)} = \left\lbrace B \in \cP_{s+\ell}: \exists B' \in \cC_j^{(\ell)} \mbox{ such that }B\cap \varphi(B') \neq \emptyset\right\rbrace.
		\end{align*}
	\end{enumerate}
\end{algorithm}

Observe that the unions
\[
C_j^{(\ell)} = \bigcup_{B \in \cC_j^{(\ell)}} B
\]
form a nested sequence in $\ell$, i.e.,
\[
C_j^{(0)}\supset C_j^{(1)} \supset \ldots \supset C_j^{(\ell)} \ldots .
\]
In fact, it is also a nested sequence in $j$, i.e.,
\[
C_0^{(\ell)}\subset C_1^{(\ell)} \ldots \subset C_j^{(\ell)} \ldots.
\]
Due to the compactness of $Q$ the continuation in Step (3) of Algorithm~\ref{alg:continuation} will
terminate after finitely many, say $J_\ell$, steps. We denote the corresponding box covering obtained by the continuation
method by
\begin{equation*}
\label{eq:zeta}
G_\ell = \bigcup_{j=0}^{J_\ell} C_j^{(\ell)} = C_{J_\ell}^{(\ell)}.
\end{equation*}

In \cite{ZDG18} we proved that increasing $\ell$ eventually leads to convergence of $C_j^{(\ell)}$ to the subsets $W_j$ and assuming that the closure of the embedded unstable manifold $\overline{W^u(p)}$ is attractive $G_\ell$ converges to $\overline{W^u(p)}$ .

\begin{proposition}[{\cite[Proposition 5]{ZDG18}}].\label{prop:convergenceCC}\quad
	\begin{enumerate}[(a)]
		\item The sets $C_j^{(\ell)}$ are coverings of $W_j$ for all $j,\ell = 0,1,\ldots$. Moreover, for fixed $j$, we have
		\[
		\bigcap_{\ell = 0}^\infty C_j^{(\ell)} = W_j.
		\]
		
		\item Suppose that $\overline{W^u(p)}$ is linearly attractive, i.e., there is a $\lambda \in (0,1)$ and a neighborhood $U\supset Q \supset \overline{W^u(p)}$ such that
		\begin{equation*}\label{eq:basin}
		\setdist{\varphi(y),\overline{W^u(p)}}\leq \lambda~\setdist{y,\overline{W^u(p)}} \quad\forall y\in U.
		\end{equation*}
		Then the box coverings obtained by Algorithm \ref{alg:continuation} converge to the closure of the embedded unstable manifold $\overline{W^u(p)}$. That is,
		\[
		\bigcap_{\ell = 0}^\infty  G_\ell=\overline{W^u(p)}.
		\]
		
	\end{enumerate}
\end{proposition}

\section{Review of Diffusion Maps}
\label{sec:Diffusion Maps}
In the last sections it was shown that combining embedding techniques with set oriented numerical methods allows the computation of one-to-one images in $\R^k$ of attractors and manifolds of infinite dimensional dynamical systems. However, the embedding can still be high dimensional, even though the box-counting dimension is low ($k>2d_{\mathrm{box}}$). Thus, the embedded set is topologically uninformative and it is hard to identify geometrical features of the underlying attractor or manifold. To highlight these important features and possibly further decrease the embedding dimension we rely on feature extraction methods such as the concept of diffusion maps~\cite{coifman2006diffusion}, whose construction we briefly review for our purposes.

Let $X=\{x_i\}_{i=1}^m\subset \R^k$ be a finite set of sample points, called \emph{anchor points}, that (coarsely) approximate the embedded attractor $A_k\subset \R^k$ or the embedded unstable manifold $W^u(p)\subset A_k$.

Suppose $k_\varepsilon:\R^k\times \R^k \to \R$ is a rotation-invariant kernel of the following form
\[
k_\varepsilon(x_1,x_2)=h\left(\frac{\norm{x_1-x_2}^2}{\varepsilon}\right).
\]
For a given $\varepsilon>0$ and $\alpha \in [0,1]$ we construct a stochastic matrix $P_{\varepsilon,\alpha}\in \R^{m\times m}$ by
\begin{align*}
&\tilde q_i=\sum_{j=1}^m k_\varepsilon(x_i,x_j), && k_\varepsilon^{(\alpha)}(i,j)=\frac{k_\varepsilon(x_i,x_j)}{\tilde q_i^\alpha \tilde q_j^\alpha},\\
&\tilde  d_i=\sum_{j=1}^m k_\varepsilon^{(\alpha)}(i,j), && \left(P_{\varepsilon,\alpha}\right)_{ij} = \tilde p(i,j)=\frac{k_\varepsilon^{(\alpha)}(i,j)}{\tilde d_i}.
\end{align*}
The choice of $\varepsilon$ and $\alpha$ will be discussed later. Observe that $P=P_{\varepsilon,\alpha}$ has a sequence of decreasing eigenvalues $\lambda_l$ and corresponding eigenvectors $\psi_l$ where $\lambda_0
=1$. Then, according to \cite{coifman2006diffusion} (or Theorem~\ref{thm:R05}) the $k$-dimensional diffusion map
\begin{align*}
\tilde R:X\to \R^k,~ x_i\mapsto y_i:=(\lambda_1\psi_1(x_i),\ldots,\lambda_k \psi_k(x_i))
\end{align*}
embeds the data into $\R^k$ (up to some relative error), where $\psi_l(x_i)$ is the $i$-th entry of the $\ell$-th eigenvector of $P$. Since this map is only defined on some data points $X=\{x_i\}_{i=1}^m$, we extend this map $\tilde R$ to a map $R:\R^k\to \R^k,~x\mapsto y$ in a natural way that is inspired by Nystr\"oms method \cite{coifman2006geometric,bengio2004out}.
For $x\in \R^k$ let
\begin{align*}
k_j(x) = k_\varepsilon(x,x_j),~q=\sum_{j=1}^m k_j \mbox{ and } k_j^{(\alpha)}=\frac{k_j}{q^\alpha \tilde q_j^\alpha}
\end{align*}
and again normalize by
\begin{align*}
d=\sum_{j=1}^m k_j^{(\alpha)} \mbox{ and } p_j=\frac{k_j^{(\alpha)}}{d}.
\end{align*}
We define the $\ell$-th entry of $y=:R(x) \in \R^k$ and $\psi_\ell(x)$, respectively, by
\begin{align}\label{eq:y}
y^{(\ell)}&:=\sum_{j=1}^m p_j \psi_\ell(x_j),\\\notag
\psi_\ell(x)&:=\frac{y^{(\ell)}}{\lambda_\ell}.
\end{align}
Note that this construction is consistent with the definition on the data set $X$, i.e., $R\vert_{X}=\tilde R$. The reason why we use this extension method is that we want to use diffusion maps not only on the given (coarse) data points but also on new data points without the costly recomputation of the whole diffusion maps. Therefore, we can easily embed trajectories of the underlying dynamical system to reveal the dynamics in diffusion coordinates or add additional data points to obtain a finer discretization.

\begin{remark}\quad
	\begin{enumerate}[(a)]
		\item In practice we use a kernel that has the form ${h(z)=c_r\exp(-z)\mathbf{1}_{z\leq r}}$ with some cutoff radius $r>0$ and constant $c_r$ such that $\int h(\norm{z}^2)\mathrm{d}z=1$ to increase the sparsity of $P$ and reduce the numerical effort. For simplicity we choose $r=\sqrt{2\varepsilon}$, to assure that interaction between data points further apart than $r$ is sufficiently small.
		\item  If we embed an out-of-sample point $x\in \R^k$ that is not in the original data set $X$, it is possible, that there is no anchor point in the $r$-ball of $x\in \R^k$ and thus $x$ will be mapped to the origin. To prevent this phenomena we have adapted the extension method. In the following we increase $r$ successively by $10\%$ for exactly those points until there are $N$ neighbors without changing~$\varepsilon$ to obtain a coefficient vector $p_j$ that has at least $N$ non-vanishing entries. 
		However, this idea is not optimal. In fact, the proposed extension method is only accurate for points within the kernel bandwidth \cite{LONG2017}.
	\end{enumerate}
\end{remark}

To find a good choice of $\varepsilon$ we rely on the observations in \cite{coifman2008graph}. They noted, when $\varepsilon$ is well tuned, the kernel localizes the data set such that
\begin{equation} \label{eq:Seps}
S(\varepsilon):=\frac{1}{m^2}\sum_{i,j}k_\varepsilon(x_i,x_j)\approx \frac{(4\pi\varepsilon)^{d_{\mathrm{int}}/2}}{\vol{(W_u)}},
\end{equation}
where $d_{\mathrm{int}}$ is the intrinsic dimension of the embedded manifold $W_u$. Therefore, $S(\varepsilon)$ should be locally well approximated by a power law $S(\varepsilon)\sim \varepsilon^a$, where 
\[
a=\frac{d(\log S)}{d(\log \varepsilon)}
\]
is the local slope at appropriate values~$\varepsilon$ for $\log S$ versus $\log \varepsilon$. Thus, in \cite{berry2016variable} $S$ was evaluated for a large range of $\varepsilon_i=2^i$ and the finite differences
\[
a_i=\frac{\log S(\varepsilon_{i+1})-\log S(\varepsilon_i)}{\log \varepsilon_{i+1}-\log \varepsilon_i}
\]
were maximized to find an ``optimal''~$\varepsilon^*$. The intrinsic dimension is then given by $d_{\mathrm{int}}=2a_{\text{max}}$ and a good choice for $\varepsilon$ would be a value near the maximizer $\varepsilon^*$ in the region of linearity. For smoother results we first find this region by analyzing the behavior of $S$ and then use a finer discretization in $\varepsilon$ inside that region to determine a suitable $\varepsilon$ and the dimension. In this process we also fix the cutoff radius as $r=\sqrt{2\max_i{\varepsilon_i}}$ to decrease the numerical effort of computing $S(\varepsilon)$.
\clearpage
Concerning the choice of $\alpha \in [0,1]$ we summarize the result of \cite{coifman2006diffusion} as follows. Suppose the data set $X$ is on an entire compact $C^{\infty}$ submanifold $\cM$ of $\R^n$ (they also discuss finite $X$ that approximates $\cM$) and is distributed\footnote{The formally correct way to state this is that the empirical measure of the data points converges weakly to a measure with density~$q$ as $m\to \infty$. The corresponding convergence results can be found in~\cite{HAvL07}.} with density $q$ and $\Delta$ is the (positive semi-definite) Laplace--Beltrami operator on $\cM$. Then $\Delta$ has eigenfunctions that verify the Neumann condition at the boundary $\partial M$ and form a Hilbert basis of $L^2(\cM,\mathrm{d}x)$. Let $E_K$ be the linear span of the first $K+1$ Neumann eigenfunction of $\Delta$.


\begin{proposition}[{\cite[Theorem~2]{coifman2006diffusion}}]
	Let
	\[
	L_{\varepsilon,\alpha} = \frac{\mathbf{I}-P_{\varepsilon,\alpha}}{\varepsilon}
	\]
	be the (discrete-time) infinitesimal generator of the Markov chain. Then for a fixed $K>0$, we have for $f\in E_K$
	\begin{align*}
	\lim_{\varepsilon\to 0} L_{\varepsilon,\alpha}f = \frac{\Delta(f q^{1-\alpha})}{q^{1-\alpha}}-\frac{\Delta(q^{1-\alpha})}{q^{1-\alpha}}f.
	\end{align*}
	In other words, the eigenfunctions of $P_{\varepsilon,\alpha}$ can be used to approximate those of the following symmetric Schr\"odinger operator:
	\[
	\Delta \phi - \frac{\Delta(q^{1-\alpha})}{q^{1-\alpha}}\phi,
	\]
	where $\phi=fq^{1-\alpha}$.
\end{proposition}

In particular, for $\alpha=1$ they found
\[
\lim_{\varepsilon\to 0} L_{\varepsilon,1} = \Delta
\]
and for any $t>0$, the Neumann heat kernel $e^{-t\Delta}$ can be approximated on $L^2(\cM)$ by $P_{\varepsilon,1}^{\frac{t}{\varepsilon}}$:
\[
\lim_{\varepsilon\to 0} P_{\varepsilon,1}^{\frac{t}{\varepsilon}}= e^{-t\Delta}.
\]
Thus for $\alpha=1$ the Markov chain converges to the Brownian motion on $\cM$. Consequently the normalization removes the influence of the density and we recover the Riemannian geometry of the data set as desired.

\section{Application of diffusion maps to embedded attractors and manifolds}
\label{sec:applications}

\subsection{Embedded manifolds of the Kuramoto--Sivashinsky equation}\label{sec:Kuramoto-Sivashinsky}
In the recent work~\cite{ZDG18} the embedded unstable manifold of $u^*=0$ of the Kuramoto--Sivashinsky equation
\begin{equation}
\label{eq:KS}
\begin{aligned}
&u_t + 4u_{xxxx}+\mu\left[ u_{xx} + \cfrac 1 2 (u_x)^2 \right] = 0, \quad 0\leq x\leq 2\pi,\\
&u(x,0) = u_0(x), \quad u(x+2\pi,t) = u(x,t)
\end{aligned}
\end{equation}
was approximated for several parameter values $\mu>0$ using Algorithm~\ref{alg:continuation}. For the construction of the observation map a POD basis $\{\zeta_1,\ldots,\zeta_S\}$ has been computed by doing a singular value decomposition on a snapshot-matrix obtained by a long-term integration of
\[
u_0=0.0001\cos(x)(1+\sin(x)).
\]
Then the observation map was given as the projection of a state $u\in Y$ onto the first $k$ POD coefficients  $\alpha_i=\langle u,\zeta_i \rangle$ for $i=1,\ldots,k$, i.e.,
\[
R(u)=(\alpha_1,\ldots,\alpha_k)^T=(\langle u, \zeta_1 \rangle,\ldots,\langle u, \zeta_k \rangle)^T,
\]
where $\langle \cdot,\cdot\rangle$ denotes the $L^2$ scalar product. For the purpose of comparing the parameter dependent manifolds in POD and diffusion coordinates we embed the manifolds with respect to the basis that is computed for $\mu=15$ if not said otherwise. However, in general the basis should be adapted to the current considered parameter.

To decrease the numerical effort we alter the original continuation algorithm. We skip step 1) and do only one continuation step but with a huge amount of test points ($100000$) and a relatively long integration time of $T=800$. However, we are not only adding those boxes that are hit after time $T$ but also the boxes that the computed embedded trajectories cross, i.e., we additionally add the boxes that contain points that are integrated at time instances $ih~,i=1,\ldots,T/h$, where $h=0.2$. 

To apply diffusion maps on the generated box covering we choose as anchor points the mid-points of $m=100000$ random boxes and approximate the optimal value $\varepsilon^*$ as described in Section~\ref{sec:Diffusion Maps} for an optimal performance of the embedding technique. Note that we under-sample the manifold  in this way and thus may underestimate the intrinsic dimensions. After computing the diffusion coordinates of the anchor points we additionally embed up to $500000$ of the remaining midpoints via the extension scheme of Nystr\"om \eqref{eq:y} to increase the density of the point cloud in diffusion coordinates.

\subsubsection{The travelling wave} For $\mu=15$ the Kuramoto--Sivashinsky equation has two stable traveling waves (limit cycles) traveling in opposite directions due to the symmetry imposed by the periodic boundary conditions. In the observation space this corresponds to two stable limit cycles that are symmetric in the first POD coefficient~$a_1$. In addition to that, a loop of unstable steady states that surrounds $u^*=0$ was found numerically by the long-term simulation for the constructing of the POD basis. Topologically, it is an entire circle due to the periodic boundary conditions in~\eqref{eq:KS}. Thus, a long-term simulation first approaches a point on this circle and then eventually converges to one of the traveling waves (the limit cycle).

\begin{figure}[H]
	\subfigure[]{\includegraphics[width = 0.49\textwidth]{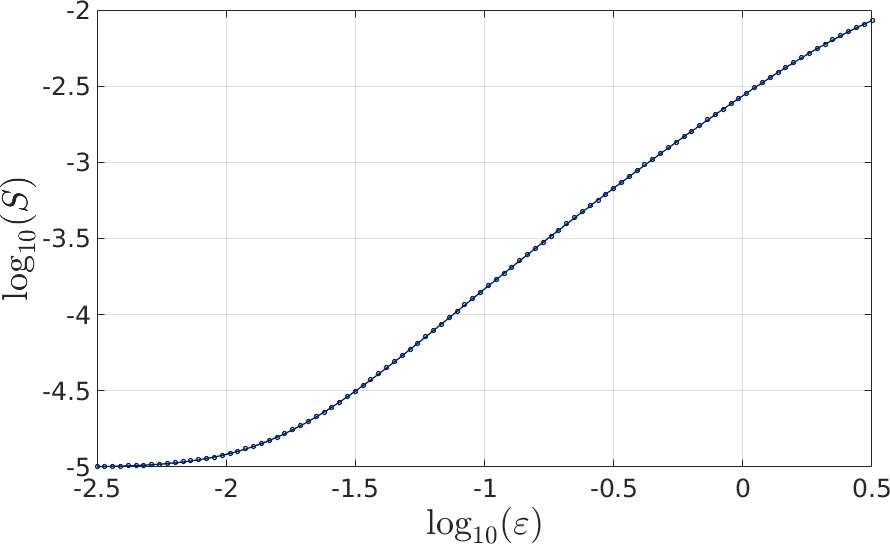}}
	\subfigure[]{\includegraphics[width = 0.49\textwidth]{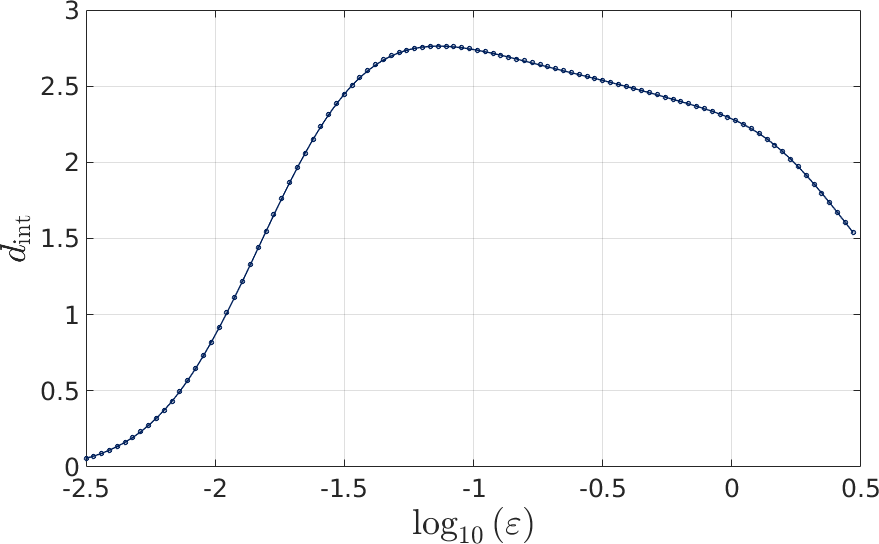}}
	\caption{(a) $\log S(\varepsilon)$ versus $\log \varepsilon$ plot, cf.~\eqref{eq:Seps}. (b) The estimated intrinsic dimension for $\mu=15$ by local slope approximation. $d_{\mathrm{int}}$ is maximized at $\varepsilon^*\approx 0.07$.}
	\label{fig:KS_15_d}
\end{figure}

\begin{figure}[H]
	\subfigure[]{\includegraphics[width = 0.49\textwidth]{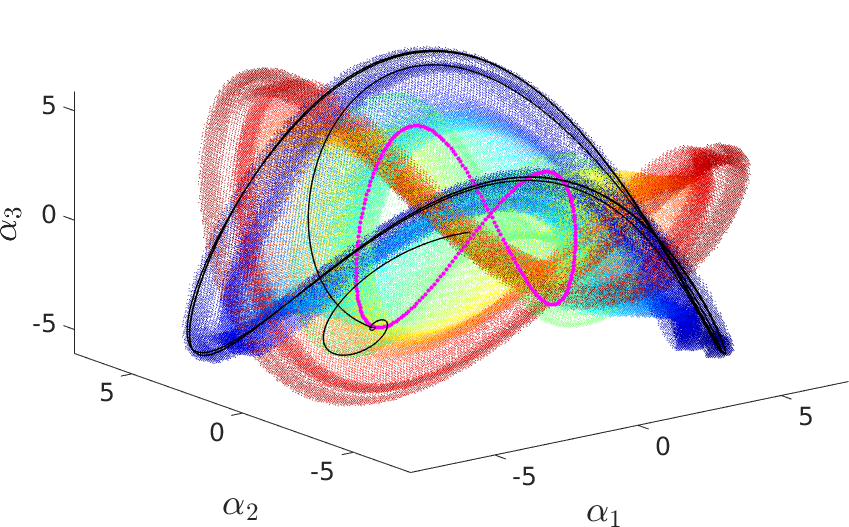}}
	\subfigure[]{\includegraphics[width = 0.49\textwidth]{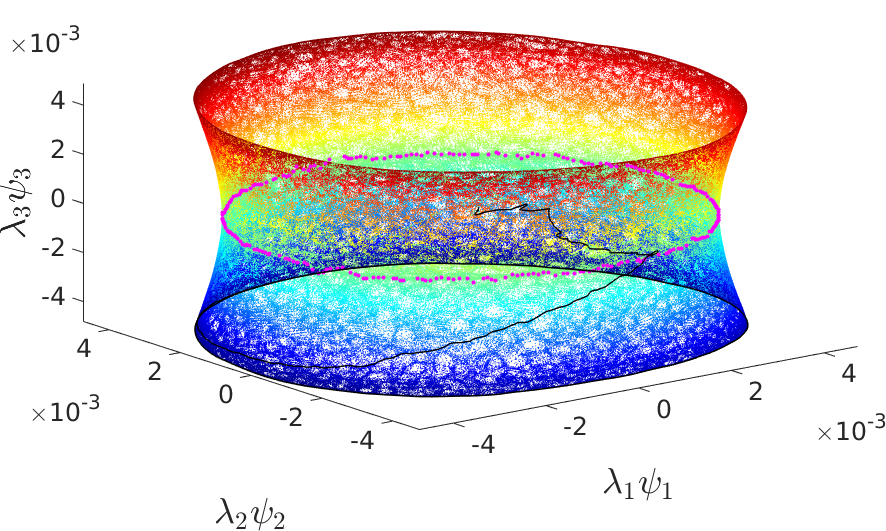}}
	\caption{POD coordinates (a) and diffusion coordinates (b) of the unstable manifold of the Kuramoto--Sivashinsky equation for $\mu=15$. The coloring is according to the third diffusion coordinate $\lambda_3\psi_3$. The embedding of the loop of unstable steady states and a trajectory starting around the unstable steady state $u^*=0$ are plotted in magenta and black, respectively.
	}
	\label{fig:KS_15}
\end{figure}

We choose an embedding dimension of $k=7$ and approximate the embedded unstable manifold at level $s=56$ with $1181433$ boxes. With the ideas used in \cite{coifman2008graph} and \cite{berry2016variable} we find $\varepsilon^*\approx 0.07$ (cf. Figure \ref{fig:KS_15_d}), where 
\[
W_u\subset[-8,8]\times[-8,8]\times[-7,7]\times[-6,6]\times[-2,2]\times[2,2]\times[0.5,0.5] \subset \R^7.
\]
Observe that our estimated dimension of at least $d_{int}\approx 2.75$ is greater than two which is caused by outer approximation. Considering the previous discussion the manifold should have a dimension of exactly two that connects $u^*=0$ with the loop of unstable steady states. However, the continuation Algorithm~\ref{alg:continuation} does not stop when that orbit is discovered.
In Figure \ref{fig:KS_15} we show the discretized embedded manifold and its diffusion coordinates. We see that the first two diffusion coordinates like the first two POD coordinates form a circular disc. But the third diffusion coordinate reveals more structure than the third POD coordinate. In fact it distinguishes between both limit cycles: $\psi_3>0$ represents convergence to the first limit cycle located at $\psi_3\approx 0.0045$, where analogously $\psi_3<0$ shows the convergence to the second limit cycle at $\psi_3\approx - 0.0045$.
In addition to that $\psi_3=0$ marks the inner part of the manifold which connects the unstable steady state $u^*=0$ with the entire orbit of unstable steady states (plotted in magenta), that lie at the boundary of the disk. We observed that the higher order coordinates are so--called higher harmonics, i.e., functions of the first three diffusion coordinates and thus not giving any additional topological information. In conclusion, the shape of the manifold can be described as a cylinder that has a disk inside it cutting it perpendicularly to its cylindrical axis.
\subsubsection{The stable heteroclinic cycle}
\begin{figure}[H]
	\subfigure[]{\includegraphics[width = 0.49\textwidth]{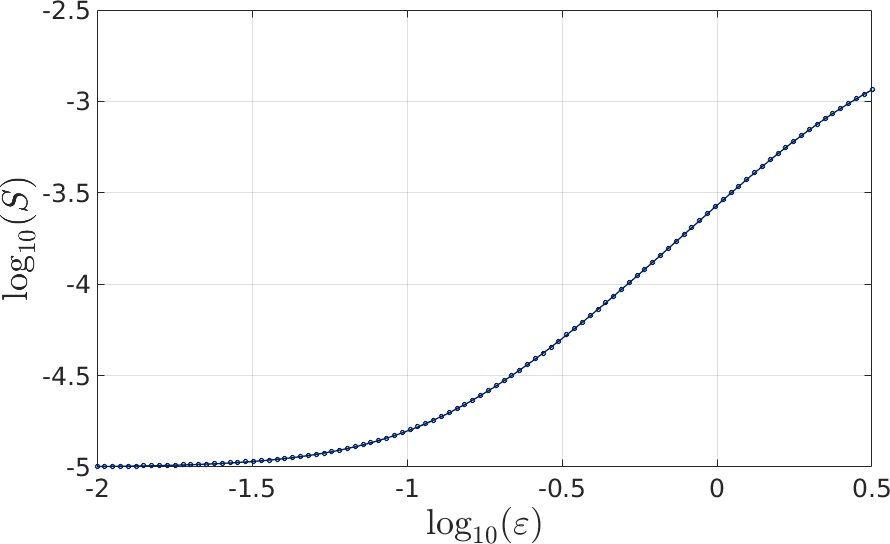}}
	\subfigure[]{\includegraphics[width = 0.49\textwidth]{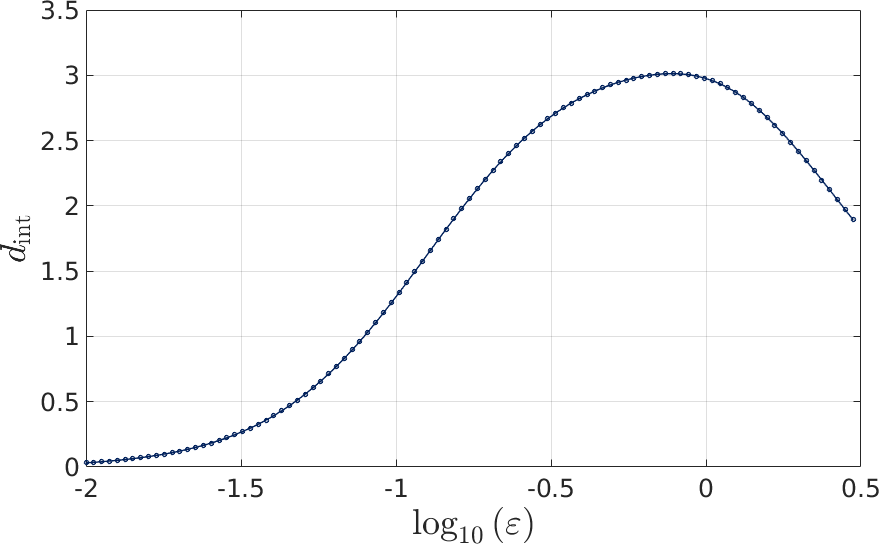}}
	\caption{(a) $\log S$ versus $\log \varepsilon$ plot. (b) The estimated intrinsic dimension for $\mu=18$ by local slope approximation. $d_{\mathrm{int}}$ is maximized at $\varepsilon^*\approx 0.78$.}
	\label{fig:KS_18_d}
\end{figure}

The long-term behavior of the Kuramoto--Sivashinsky equation for $\mu=18$ is described by a pulsation between two unstable states, that are $\pi/2$-translations of each other. Moreover, the transients stay close to one of the states for a relative long time until it pulses back to the other state. Due to the boundary conditions translations of this states are also unstable states and thus different pulsations resulting from different initial conditions give trajectories between different unstable states on that loop. In fact, they are rotations of one another about the origin. The $1180913$ boxes covering the embedded unstable manifold $W_u\subset \R^7$ generated by the continuation method at level $s=35$ approximate an at least $3$-dimensional set (cf. Figure~\ref{fig:KS_18_d}), where
\[
W_u\subset [-8,8]\times[-8,8]\times[-6.5,6.5]\times[-6,6]\times[-2,2]\times[-2,2]\times[-0.5,0.5] \subset \R^7.
\]

\begin{figure}[H]
	
	\subfigure[]{\includegraphics[width = 0.49\textwidth]{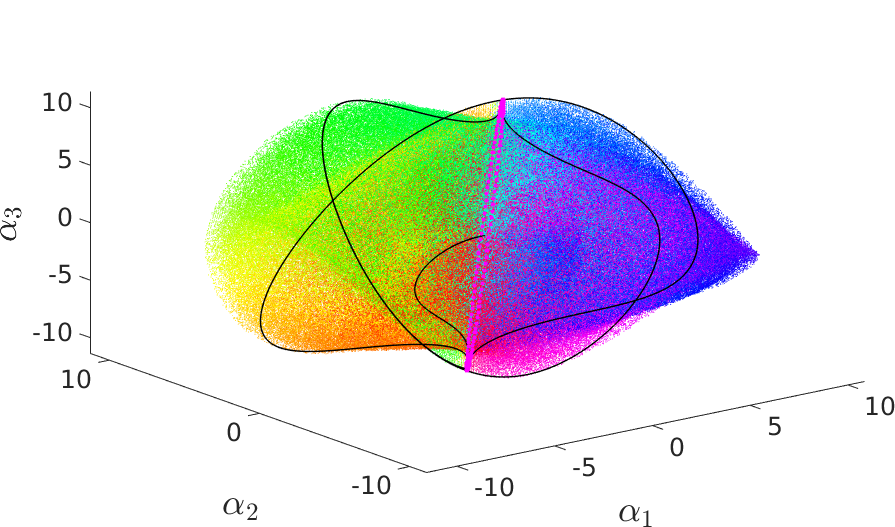}}
	\subfigure[]{\includegraphics[width = 0.49\textwidth]{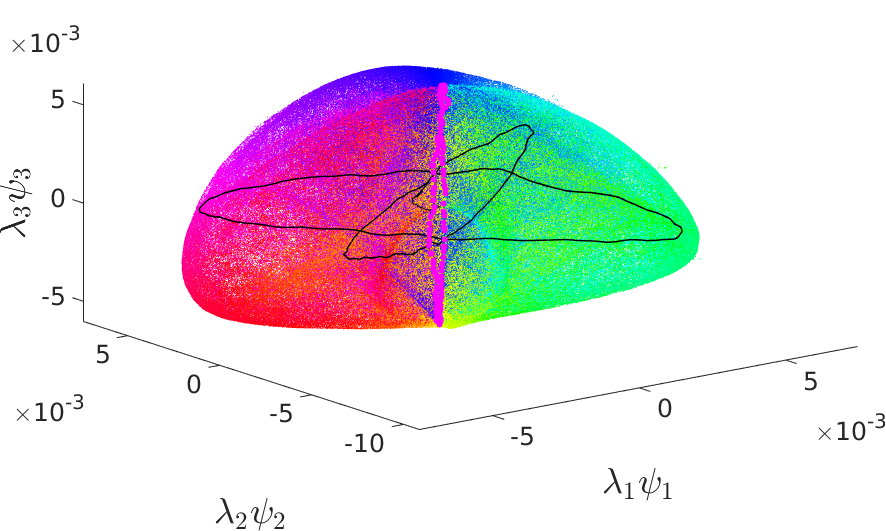}}
	\subfigure[]{\includegraphics[width = 0.49\textwidth]{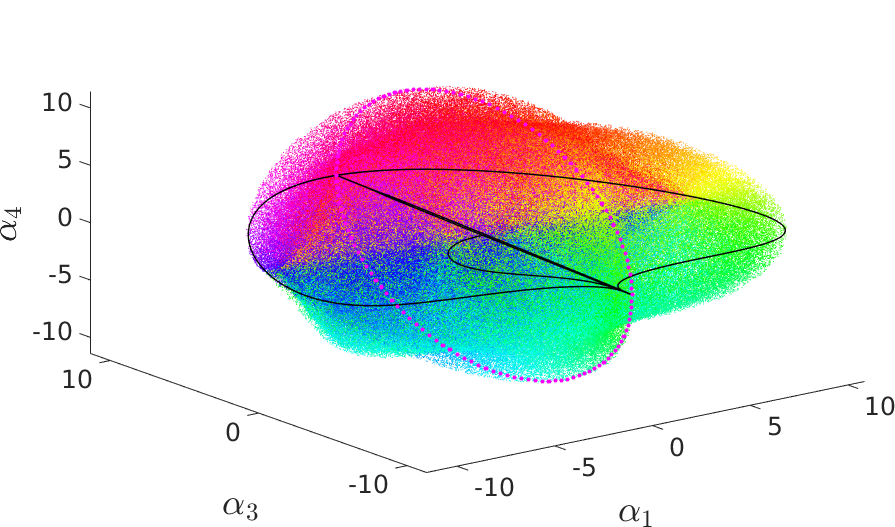}}
	\subfigure[]{\includegraphics[width = 0.49\textwidth]{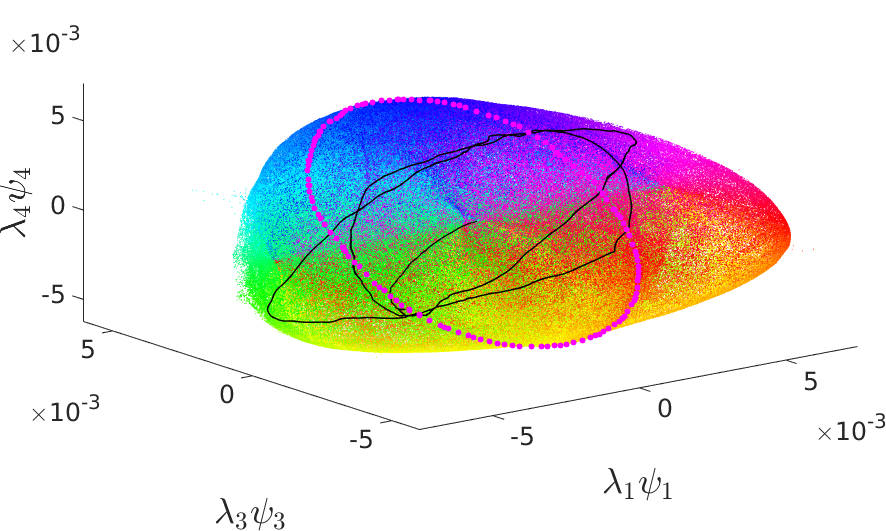}}
	\caption{POD coordinates (left column (a,c)) and diffusion coordinates (right column (b,d)) of the unstable manifold of the Kuramoto--Sivashinsky equation for $\mu=18$. The coloring is corresponds to the phase in the $\psi_1-\psi_2$ and $\psi_3-\psi_4$, respectively. A long-term simulation representing the pulsation behavior is also plotted in black and the loop of unstable states is shown magenta.}
	\label{fig:KS_18}
\end{figure}

We compute an optimal value for $\varepsilon$ of approximately $0.78$ (cf. Figure~\ref{fig:KS_18_d}) and show the corresponding embedding of the data set with respect to different diffusion coordinates in Figure~\ref{fig:KS_18}. The manifold for $\mu=18$ strongly changes its shape in POD coordinates (as expected) and also in diffusion coordinates compared to the cylindrical shape for $\mu=15$. We see, that the trajectory of the long-term simulation pulse between two states and its transients bound the embedded manifold in POD and diffusion coordinates. Also observe that the loop of unstable states is almost a straight line in the projection on $(\alpha_1,\alpha_2,\alpha_3)$ and $(\psi_1,\psi_2,\psi_3)$ (Figure~\ref{fig:KS_18} (a,b)), respectively, but is clearly visible in the $\alpha_3-\alpha_4$ and $\psi_3-\psi_4$ plane (Figure~\ref{fig:KS_18} (c,d)). Hence, we conclude that the dimension of the manifold is at larger than three and we under-sampled the manifold (cf. Figure~\ref{fig:KS_18_d}). Indeed, we find a dimension of four (see Figure~\ref{fig:KS_bif_d} (b)).


\subsection{Bifurcation analysis} 
\begin{figure}[H]
	\subfigure[]{\includegraphics[width=0.49\textwidth]{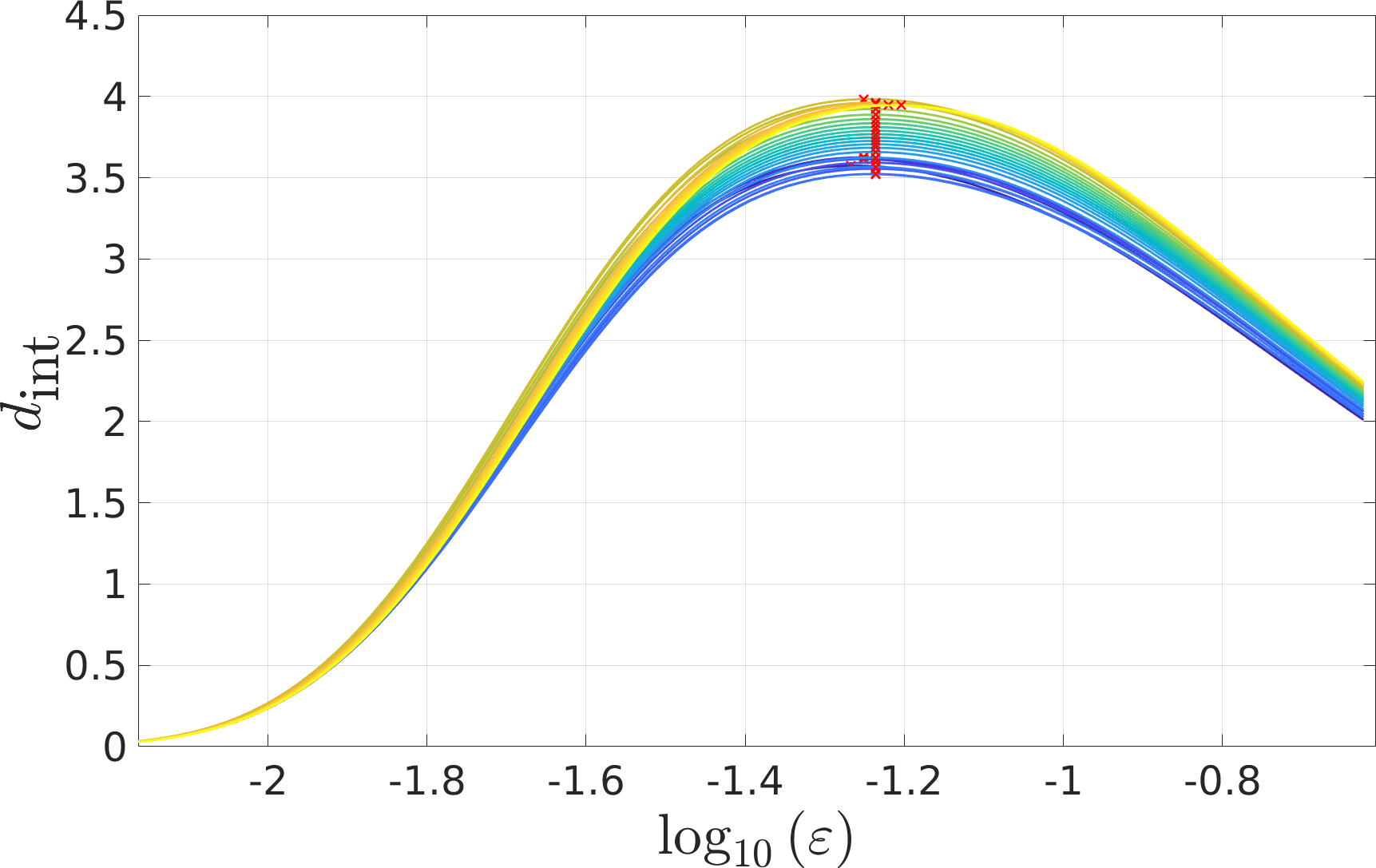}}
	\subfigure[]{\includegraphics[width=0.49\textwidth]{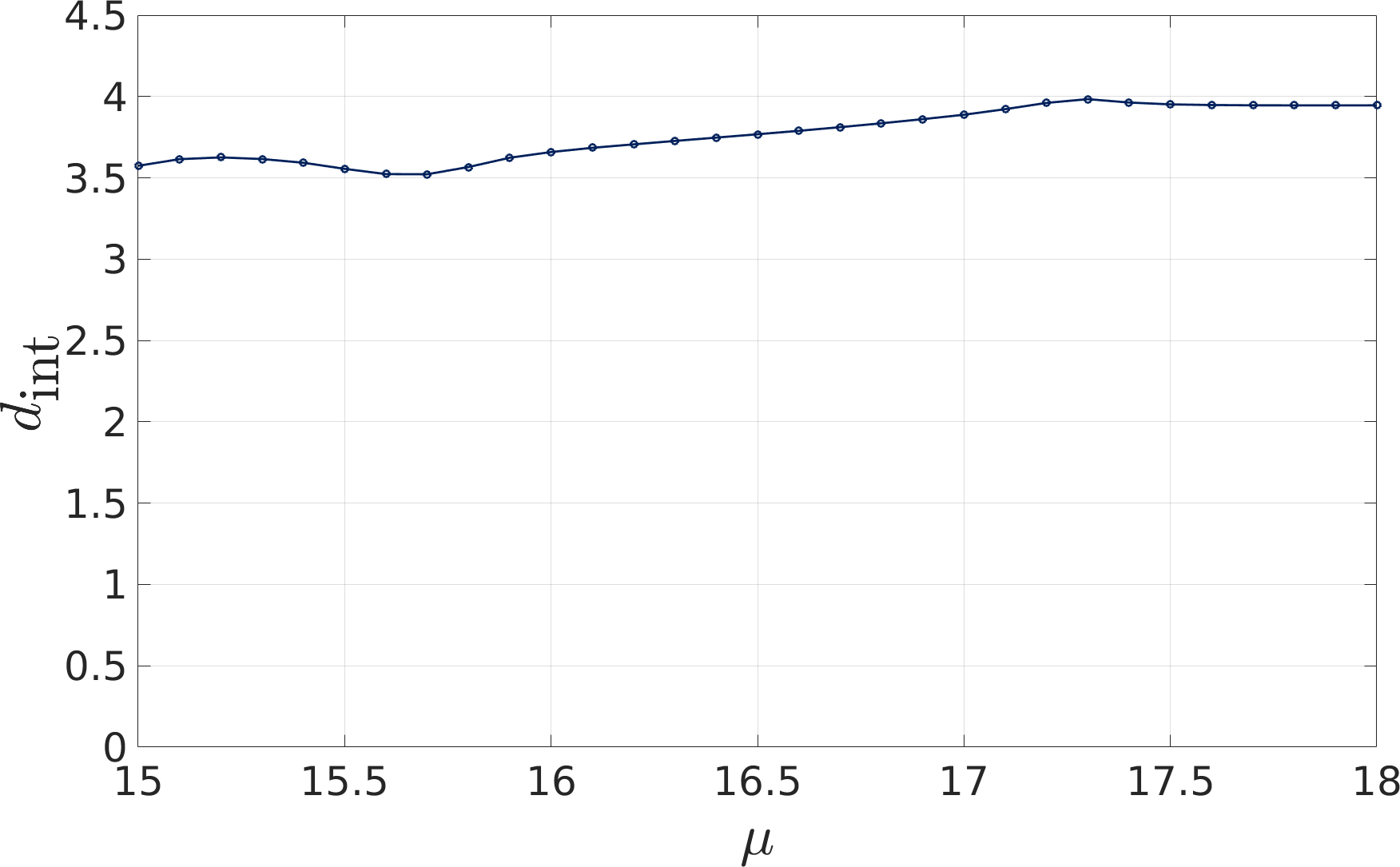}}
	\caption{Development of the intrinsic dimension of the coarse embedded manifold of the Kuramoto--Sivashinsky equation for $\mu\in[15, 18]$. (a) Local slope approximations of the intrinsic dimension. The coloring is from $\mu=15$ (blue) to $\mu=18$ (yellow). Red crosses indicate the local maxima. (b) Intrinsic dimension as a function of $\mu\in[15, 18]$.}
	\label{fig:KS_bif_d}
\end{figure}

Previous research \cite{hyman1986order} and our observation show that the unstable manifold strongly changes its structure depending on the parameter $\mu$. To further investigate this behavior we will analyze how the cylindrical shape that is revealed in diffusion coordinates for $\mu=15$ changes by increasing the parameter. Thus, our focus lies in following these three coordinates and we neglect new appearing diffusion coordinates with larger eigenvalues. In this work we select the appropriate eigenvectors by hand, but this should be algorithmically improved by a path following method in future research.
Instead of drawing some random mid-points as anchor points, we choose a coarser discretization by considering all midpoints of the approximation at level $s=49$. Again, this reduces the number of anchor points but in addition to that deals with the problem of possibly sampling the manifold poorly. Another advantage is that the anchor points lie on a grid and we can identify an uniformly optimal $\varepsilon\approx 0.06$ for all  $\mu$ (cf. Figure~\ref{fig:KS_bif_d}). The estimated intrinsic dimension $d_\mathrm{int}$ is larger than the previously found dimension which is due to the fact, that the manifold is not under-sampled like in the previous section.



\begin{figure}[H]
	\subfigure{\includegraphics[width=0.24\textwidth]{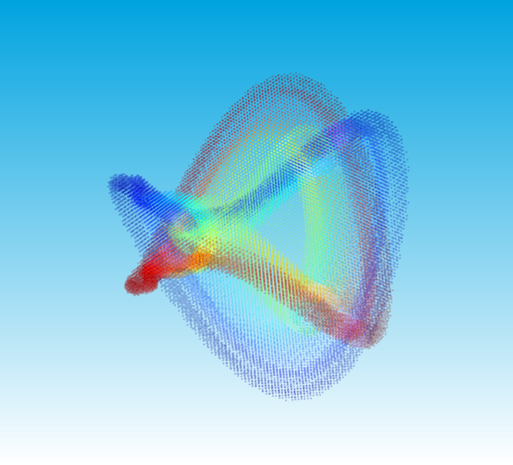}}
	\subfigure{\includegraphics[width=0.24\textwidth]{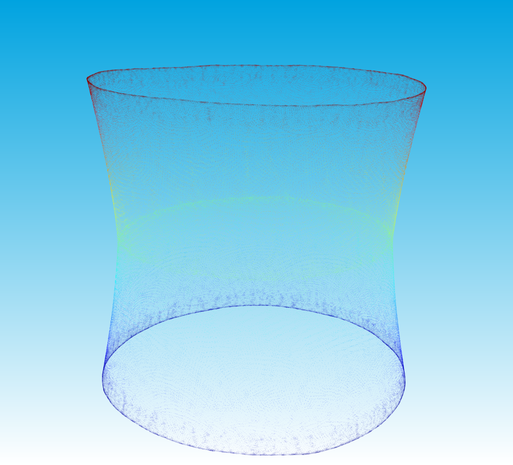}}
	\subfigure{\includegraphics[width=0.24\textwidth]{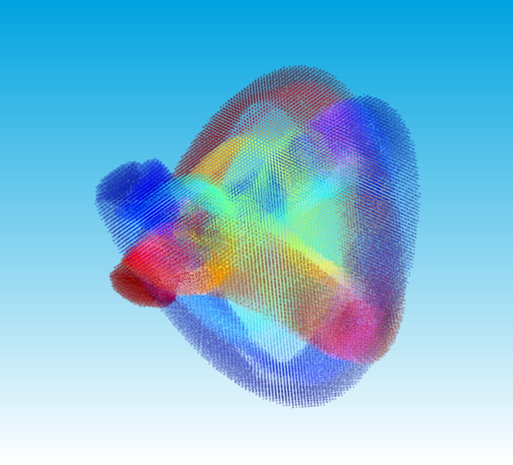}}
	\subfigure{\includegraphics[width=0.24\textwidth]{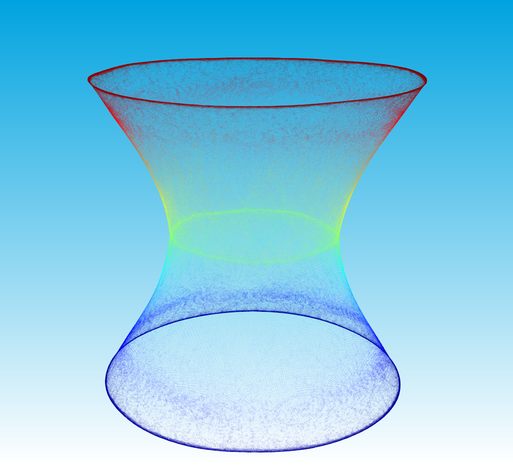}}
	\subfigure{\includegraphics[width=0.24\textwidth]{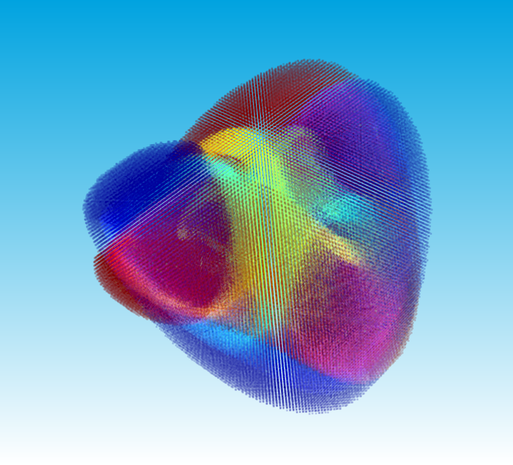}}
	\subfigure{\includegraphics[width=0.24\textwidth]{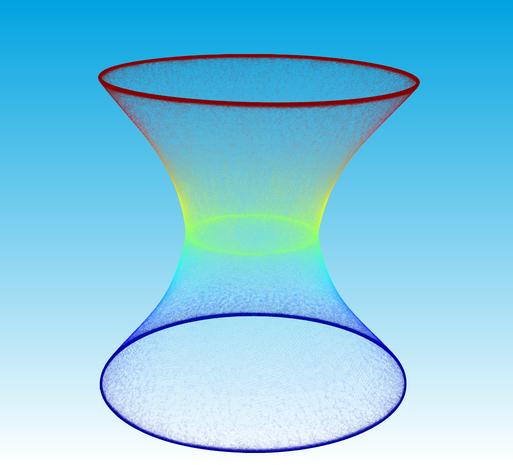}}
	\subfigure{\includegraphics[width=0.24\textwidth]{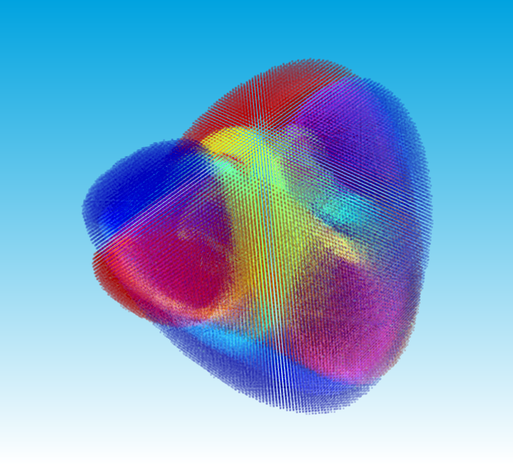}}
	\subfigure{\includegraphics[width=0.24\textwidth]{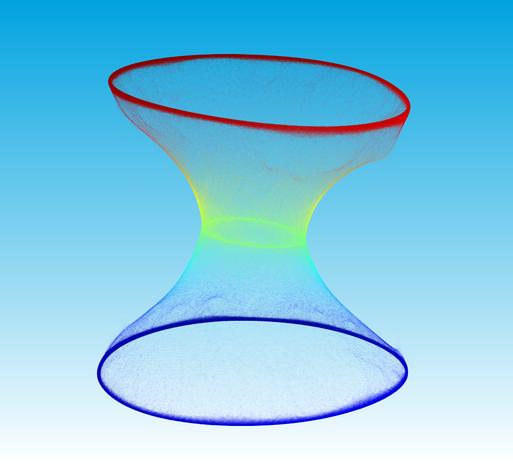}}
	\subfigure{\includegraphics[width=0.24\textwidth]{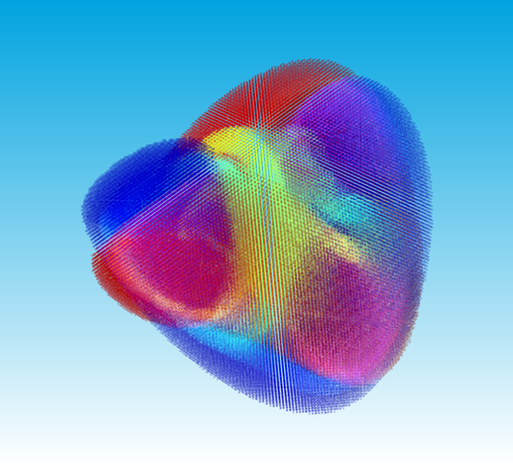}}
	\subfigure{\includegraphics[width=0.24\textwidth]{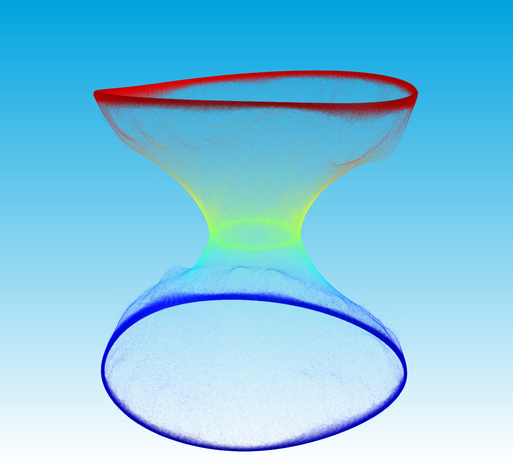}}
	\subfigure{\includegraphics[width=0.24\textwidth]{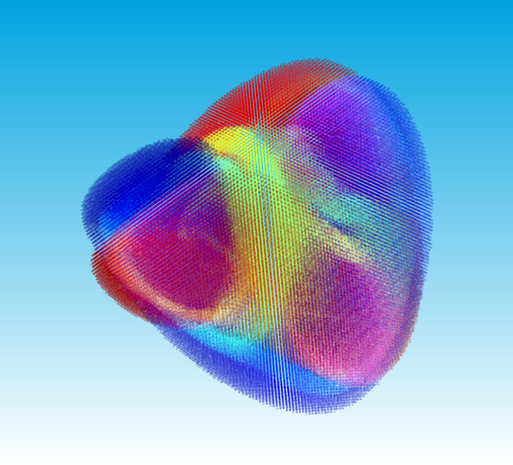}}
	\subfigure{\includegraphics[width=0.24\textwidth]{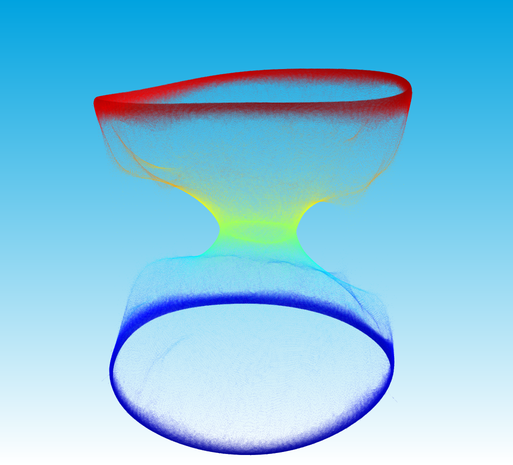}}
	\subfigure{\includegraphics[width=0.24\textwidth]{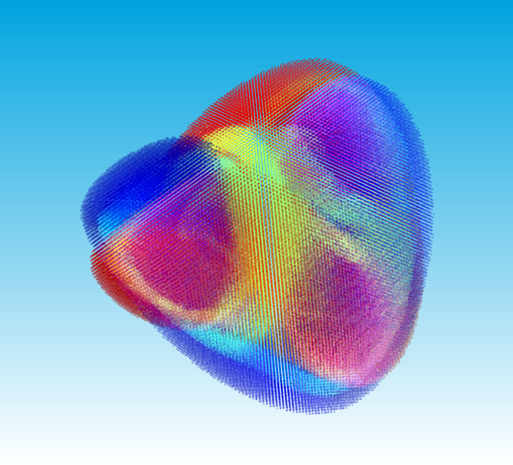}}
	\subfigure{\includegraphics[width=0.24\textwidth]{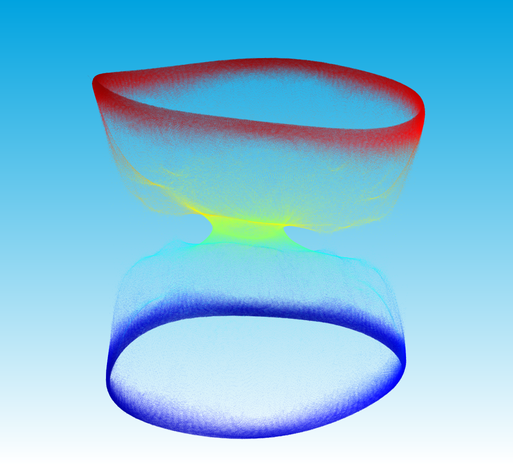}}
	\subfigure{\includegraphics[width=0.24\textwidth]{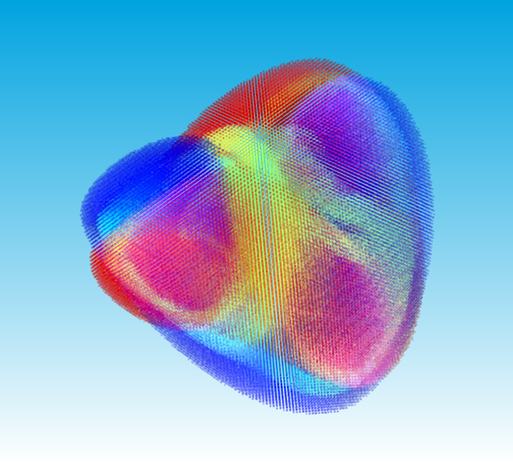}}
	\subfigure{\includegraphics[width=0.24\textwidth]{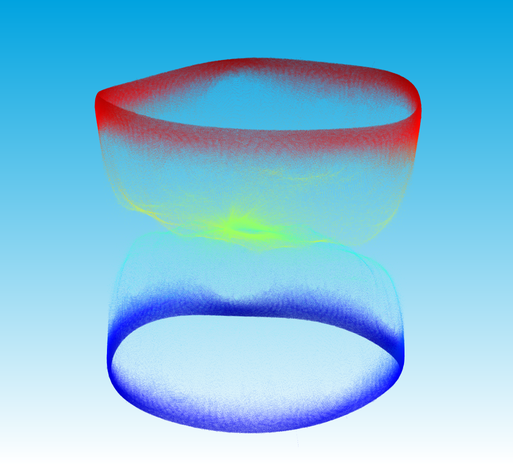}}
	\subfigure{\includegraphics[width=0.24\textwidth]{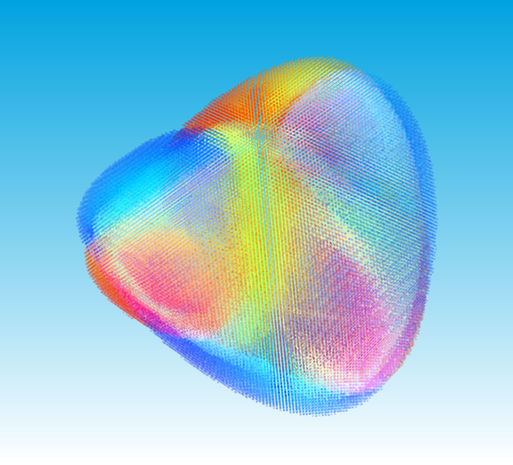}}
	\subfigure{\includegraphics[width=0.24\textwidth]{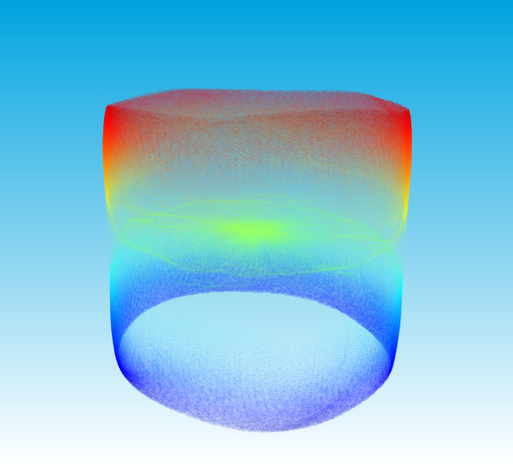}}
	\subfigure{\includegraphics[width=0.24\textwidth]{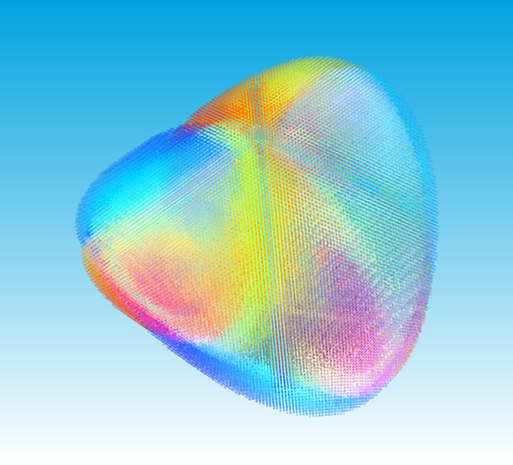}}
	\subfigure{\includegraphics[width=0.24\textwidth]{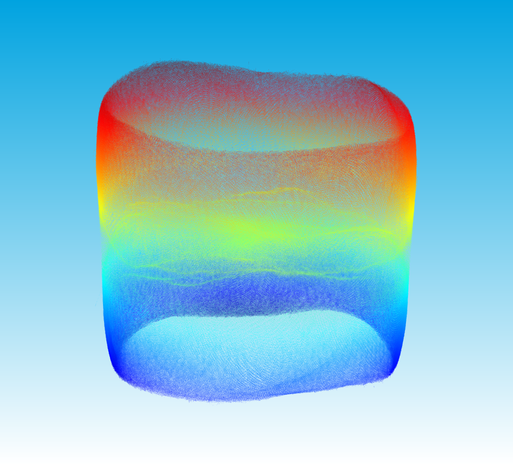}}
	\caption{Embedded manifolds of the Kuramoto--Sivashinsky equation for parameter values ${\mu\in\{15.00,15.50,16.0,16.04,16.06,16.08,16.10,16.12,16.3,16.4\}}$ in POD and diffusion coordinates (from left to right). The coloring is according to the chosen third (vertical) diffusion coordinate.
	}
	\label{fig:KS_bif}
\end{figure}

In Figure \ref{fig:KS_bif} we illustrate the changing geometry of the manifold for increasing~$\mu$. The cylindrical shape in diffusion coordinates deforms such that the circle that corresponds to the loop of unstable state shrinks together and appears to eventually bifurcate to one point in diffusion coordinates. In POD coordinates the embedded manifold becomes thicker and one quickly cannot identify the limit cycles by eye anymore. However, representing the object with by diffusion coordinates still reveals them.

As we mentioned above, while increasing $\mu$ from 15 to 18 the manifold bifurcates from a two-dimensional into a higher-than-three dimensional set. Thereby one loop of hyperbolic steady states vanishes (the pinching of the cylinder in Figure~\ref{fig:KS_bif}), and another arises (the magenta loop in Figure~\ref{fig:KS_18}). This new loop is connected by two one-parameter families of heteroclinic orbits, thus the heteroclinic orbits build two tori that intersect in one loop. A trajectory starting sufficiently close to the fixed point $p$ on the loop transitions close to, say, another fixed point $q$ on the loop by moving along one torus, then transitions back close to $p$ by moving along the other torus.

Numerical simulations show that the limit cycles (traveling waves) stay stable up to $\mu\approx 17.1$, but for $\mu\geq16$ the heteroclinic pulsation present for $\mu=18$ is a transient in the long-term behavior. Afterwards, for $\mu\geq 17.1$ the pulsation becomes dominant and convergence to a traveling wave does not occur. In future work we would like to understand whether the limit cycles bifurcate into the loop of heteroclinic points. For this we need to overcome the challenge of sufficient (and sufficiently uniform) sampling of the manifold for $\mu \geq 16.5$, that currently poses a computational bottleneck. In addition to that, as already mentioned, the selection of the correct eigenvectors has to be improved.



\subsection{The Oseberg transition}
\begin{figure}[H]
	\subfigure[]{\includegraphics[width = 0.49\textwidth]{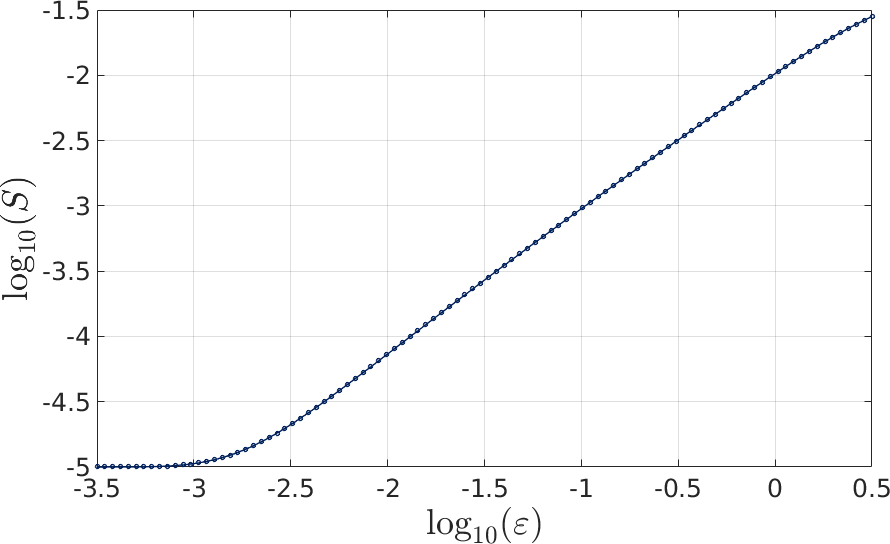}}
	\subfigure[]{\includegraphics[width = 0.49\textwidth]{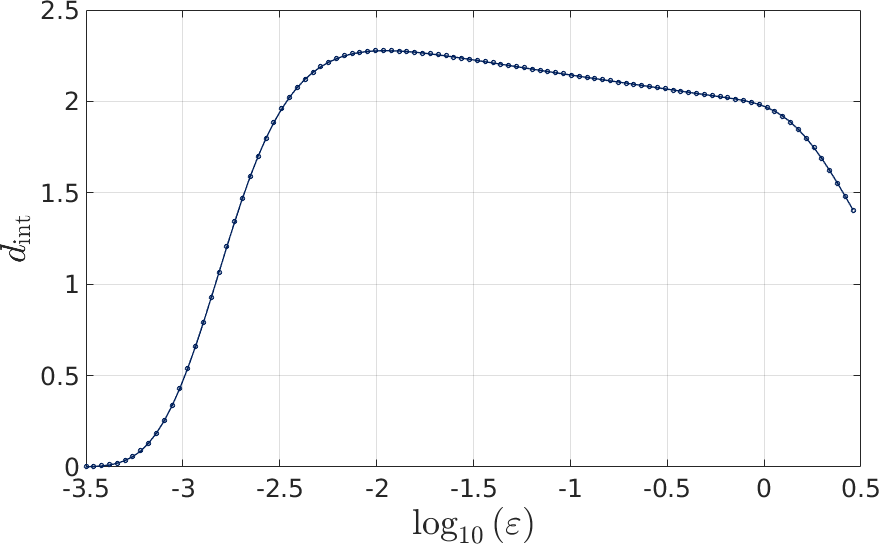}}
	\caption{(a) $\log S$ versus $\log \varepsilon$ plot. (b) The estimated intrinsic dimension for $\mu=32$ by local slope approximation. $d_{\mathrm{int}}$ is maximized at $\varepsilon^*\approx 0.01$, but we choose $\varepsilon=0.5$ to guarantee eigenvalue convergence.}
	\label{fig:KS_32_d}
\end{figure}

Finally, we consider $\mu=32$ -- the so called Oseberg transition (see~\cite{JJK01}), where the chosen initial condition near the unstable steady state $u^*=0$ is first attracted to an unstable so--called bimodal steady state, and afterwards accumulates on a limit cycle as~$t\to \infty$. Since the POD basis for $\mu=15$ is not appropriate for this parameter anymore, we adapted the basis to $\mu=32$ and observed the first $k=5$ POD coordinates with respect to that basis. We approximate the embedded manifold $W_u$ at level $s=45$ with $181443$ boxes such that
\[
W_u \subset [-13,11]\times[-9,9]\times[-6,6]\times[-2,7]\times[-2,2].
\]
To apply diffusion maps we choose $\varepsilon=0.5$ (cf. Figure~\ref{fig:KS_32_d}) since for the estimated optimal $\varepsilon^* \approx 0.01$ the convergence of eigenvalues of the diffusion matrix $P$ fails -- contrary to expectations.

\begin{figure}[H]
	\subfigure[]{\includegraphics[width = 0.49\textwidth]{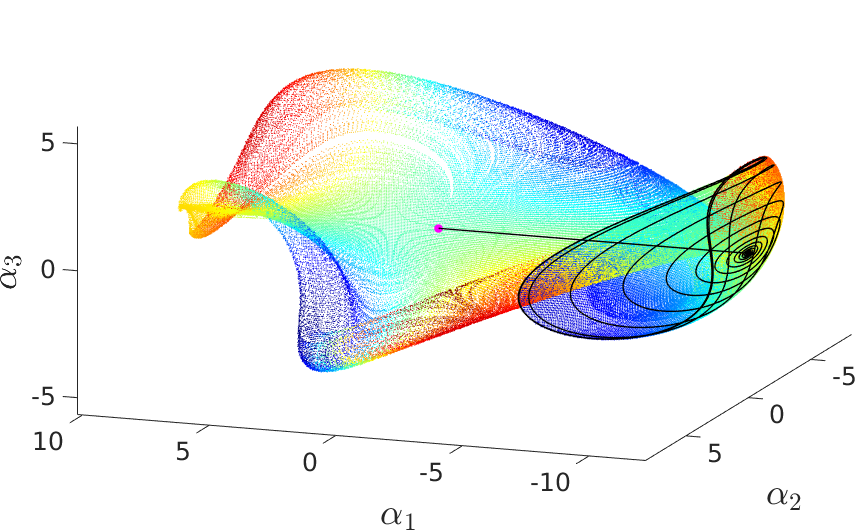}}
	\subfigure[]{\includegraphics[width = 0.49\textwidth]{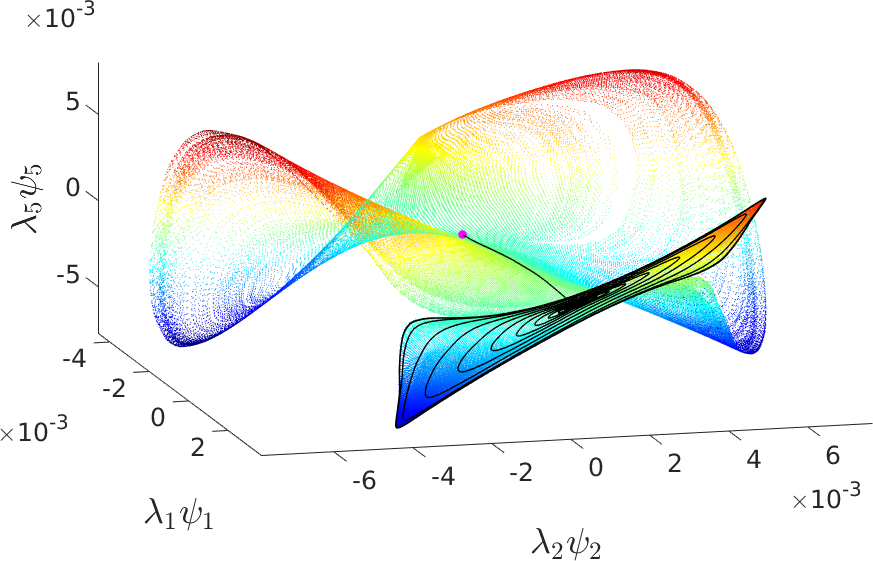}}
	\caption{POD coordinates (a) and diffusion coordinates (b) of the unstable manifold of the Kuramoto--Sivashinsky equation for $\mu=32$. The coloring is according to the diffusion coordinate $\lambda_5\psi_5$. The embedding of the unstable steady state $u^*=0$ and the orbit of unstable steady states are plotted in magenta and black, respectively.}
	\label{fig:KS32}
\end{figure}

Figure~\ref{fig:KS32} shows, how the ``jellyfish'' seen in POD coordinates is unraveled in diffusion coordinates. The corresponding long-term simulation for the computation of the POD basis is also shown in black. Observe, that we skip the third and forth diffusion coordinates since they are higher harmonics of the first and second coordinate, i.e., they are functions of the first and second diffusion coordinate and thus do not contain additional information.

\subsection{The embbeded attractor of the Mackey--Glass equation}\label{ssec:Mackey-Glass}
\begin{figure}[H]
	\subfigure[]{\includegraphics[width = 0.49\textwidth]{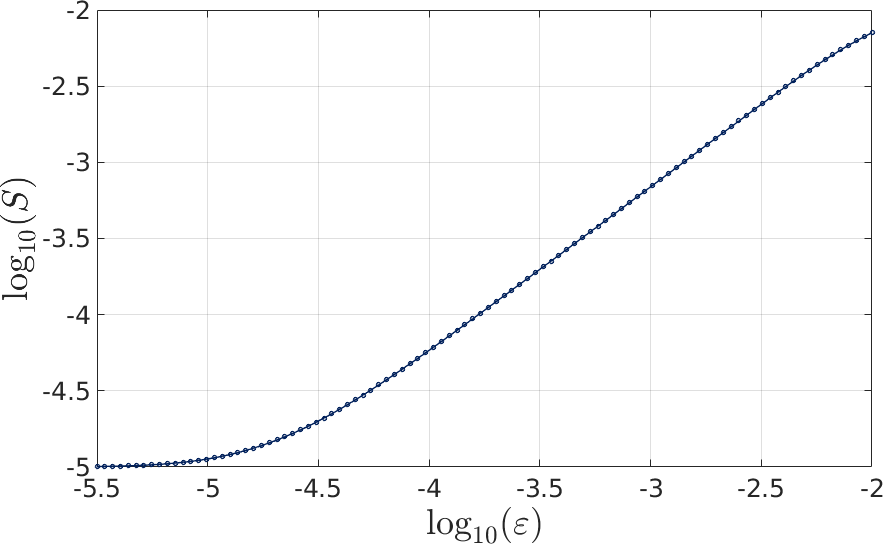}}
	\subfigure[]{\includegraphics[width = 0.49\textwidth]{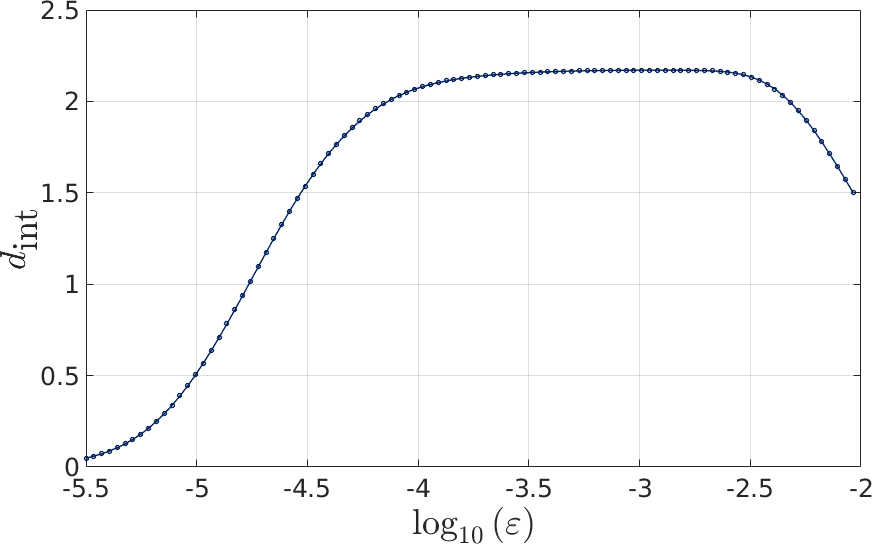}}
	\caption{(a) $\log S$ versus $\log \varepsilon$ plot. (b) The estimated intrinsic dimension of the Mackey--Glass attractor by local slope approximation. $d_{\mathrm{int}}$ is maximized at $\varepsilon^*\approx 0.0012$.}
	\label{fig:MG_d}
\end{figure}

Finally, we apply diffusion maps on the embbeded attractor of a delay differential equation with constant delay. We consider the delay differential equation introduced by Mackey and Glass in 1977~\cite{mackey1977} defined by
\begin{equation*}\label{eq:mg}
\dot u(t) = \beta \frac{u(t-\tau)}{1+u(t-\tau)^\eta}-\gamma u(t),
\end{equation*}
where we choose $\beta = 2, \gamma = 1$, $\eta = 9.65$, and $\tau = 2$. A natural observation map is given by delay coordinates
\[
R(u)=(u(-\tau),\Phi(u)(-\tau),\ldots,\Phi^{k-1}(u)(-\tau))^T.
\]
For the Mackey--Glass equation $k=7$ delays were used, i.e.
\[
R(u)=(u(-\tau),u(-\tau+\frac{\tau}{k-1}),u(-\tau+\frac{2\tau}{k-1}),\ldots,u(0))^T,
\]
to construct the core dynamical system. Then Algorithm~\ref{alg:subdivision} generated a cover of the embedded attractor with $5023208$ boxes at level $63$. Again we sample $m=100000$ random mid-points as anchor points and compute an optimal $\varepsilon^*\approx 0.0012$ (cf. Figure \ref{fig:MG_d}), where $R(u)\in [0,1.5]^7$.

\begin{figure}[H]
	\subfigure[]{\includegraphics[width = 0.49\textwidth]{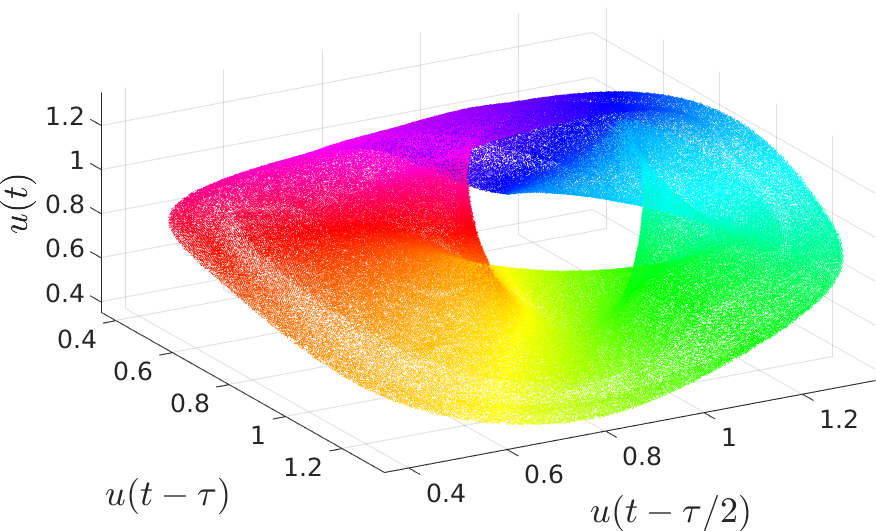}}
	\subfigure[]{\includegraphics[width = 0.49\textwidth]{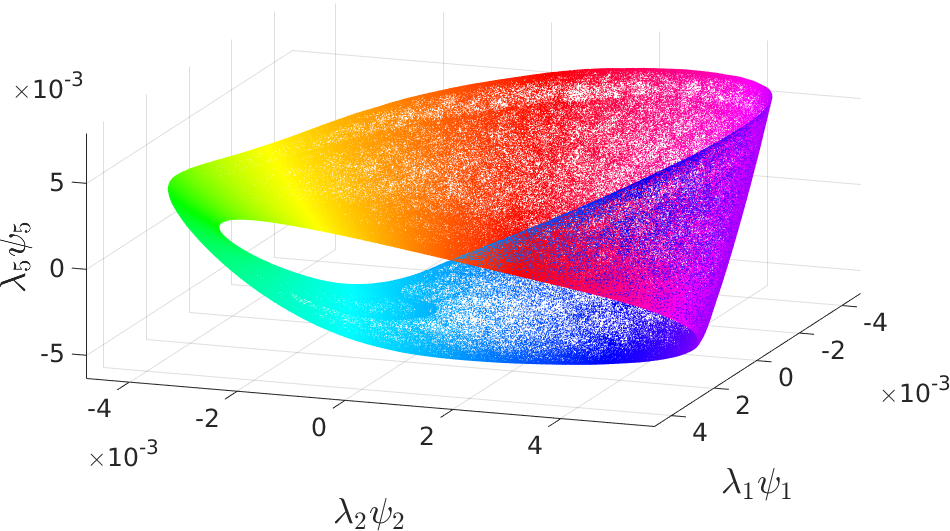}}
	\caption{Delay coordinates (a) and diffusion coordinates (b) of the attractor for the Mackey--Glass equation. The coloring is with respect to the phase in the $\psi_1-\psi_2$ plane.}
	\label{fig:MG}
\end{figure}

To increase the density of our coarse data set we again additionally embed $500000$ points via the Nystr\"om method \ref{eq:y}. The corresponding delay and diffusion coordinates are shown in Figure~\ref{fig:MG}, where the chosen diffusion coordinates reveal a Moebius strip like structure in diffusion coordinates, which is not directly clear in delay coordinates. By coloring the attractor with respect to the angle in the  $\psi_1-\psi_2$ plane we can see the phase along the strip.

\section{Conclusion}\label{sec:conclusion}

In this work we identified intrinsic coordinates of finite dimensional invariant sets of infinite dimensional systems. To this end, we first approximated these sets with set-oriented methods using observations such as POD and delay coordinates. Afterwards, we applied diffusion maps on the generated data to learn the intrinsic geometry. In future research we aim to approximate the invariant set in diffusion coordinates right away such that we construct the core dynamical system with diffusion maps as the observation map~$R$. To implement the core dynamical system numerically the extension method for out-of-sample points has to be improved to smooth the embedding and the inverse $E$ has to be numerically realized (cf.\ Section~\ref{sec:comp_emb_att}), i.e., the diffusion map embedding has to be reversed. In particular, for given $y\in \R^k$ a point $x\in \R^N$ can be computed such that $R(x)=y$ at least approximately, where $R$ is the diffusion map. Then, to combine this specific realization of the core dynamical system with set-oriented approximation techniques one has to deal with the problem of finding an initial set of anchor points and generating an initial diffusion maps embedding. For chaotic systems a long-term simulation of the system can be used, but for higher dimensions the ``uniformity'' of samples of the set will play a role.

Furthermore, it is interesting to learn the dynamics in diffusion coordinates and eventually find a reduced model (topological or in form of equations) in those coordinates~\cite{brunton2016discovering}. For instance, diffusion maps suggests that cylinder coordinates suit very well for the dynamics on the unstable manifold of the Kuramoto--Sivashinsky equation for $\mu=15$ and one might be able to construct an ordinary differential equation that describes the dynamics on the manifold.

\section*{Acknowledgments} This work is supported by the Priority Programme SPP 1881 Turbulent Superstructures of the Deutsche Forschungsgemeinschaft.

\bibliographystyle{myalpha}

\bibliography{CDS_meets_DMAPS}
\end{document}